\documentclass[11pt,twoside,leqno]{aomamlt2e}
\pageno{1}
\received{August 28, 2006}
\revised{January 27, 2009}

\renewcommand{\theequation}{\thesubsection.\arabic{equation}}

\theoremstyle{definition}

\theoremstyle{remark}

\theoremstyle{coroll}

\theoremstyle{conject}

\newcommand{\nc}{\newcommand}

\nc{\fa}{\mathfrak a} \nc{\fg}{\mathfrak g} \nc{\fk}{\mathfrak k}
\nc{\fh}{\mathfrak h} \nc{\fm}{\mathfrak m} \nc{\fn}{\mathfrak n}
\nc{\fA}{\mathfrak A} \nc{\fC}{\mathfrak C} \nc{\fI}{\mathfrak I}
\nc{\fL}{\mathfrak L} \nc{\fS}{\mathfrak S}

\nc{\nen}{\newenvironment} \nc{\ol}{\overline}
\nc{\ul}{\underline} \nc{\lra}{\longrightarrow}
\nc{\lla}{\longleftarrow} \nc{\Lra}{\Longrightarrow}
\nc{\Lla}{\Longleftarrow} \nc{\Llra}{\Longleftrightarrow}
\nc{\hra}{\hookrightarrow} \nc{\iso}{\overset{\sim}{\lra}}
\nc{\Hom}{\mathrm{Hom}} \nc{\Mor}{\mathrm{Mor}}

\nc{\notebox}[1]{\noindent\fbox{\parbox{12.5cm}{\sf #1}}\\[8pt]}
\nc{\Thm}[1]{Theorem~\ref{#1}} \nc{\Prop}[1]{Proposition~\ref{#1}}
\nc{\Lem}[1]{Lemma~\ref{#1}} \nc{\Cor}[1]{Corollary~\ref{#1}}
\nc{\Conj}[1]{Conjecture~\ref{#1}} \nc{\Claim}[1]{Claim~\ref{#1}}
\nc{\Defn}[1]{ Definition~\ref{#1}} \nc{\Exa}[1]{Example~\ref{#1}}
\nc{\Rem}[1]{Remark~\ref{#1}} \nc{\Note}[1]{Note~\ref{#1}}

\nen{thm}[1]{\label{#1}{\bf Theorem.} \em}{}
\nen{thma}[1]{\label{#1}{\bf Theorem A.} \em}{}
\nen{thmb}[1]{\label{#1}{\bf Theorem B.} \em}{}
\nen{prop}[1]{\label{#1}{\bf Proposition.} \em}{}
\nen{lem}[1]{\label{#1}{\bf Lemma.} \em}{}
\nen{klem}[1]{\label{#1}{\bf Key Lemma.} \em}{}
\nen{cor}[1]{\label{#1}{\bf Corollary.} \em}{}
\nen{cora}[1]{\label{#1}{\bf Corollary A.} \em}{}
\nen{corb}[1]{\label{#1}{\bf Corollary B.} \em}{}
\nen{conj}[1]{\label{#1}{\bf Conjecture.} \em}{}
\nen{claim}[1]{\label{#1}{\bf Claim.} \em}{}
 \nen{claima}[1]{\label{#1}{\bf Claim A.} \em}{}
\nen{claimb}[1]{\label{#1}{\bf Claim B.} \em}{}

\nen{defn}[1]{\label{#1}{\bf Definition.} }{}
\nen{exa}[1]{\label{#1}{\bf Example.}}{}

\nen{rem}[1]{\label{#1}{\em Remark.} }{}
\nen{exer}[1]{\label{#1}{\em Exercise.} }{}
\nen{pf}{\Proof}{\Endproof}
\nen{rems}[1]{\label{#1}{\em Remarks.} }{}

 \nc{\br}{\mathbb R}
 \nc{\bz}{\mathbb Z}
 \nc{\bc}{\mathbb C}
 \nc{\bn}{\mathbb N}
 \nc{\G}{\Gamma}
 \nc{\sm}{\setminus}
 \nc{\sub}{\subset}
 \nc{\lm}{\lambda}
 \nc{\al}{\alpha}
 \nc{\bt}{\beta}
 \nc{\kap}{\kappa}
 \nc{\om}{\omega}
 \nc{\dl}{\delta}
 \nc{\g}{\gamma}
 \nc{\Dl}{\Delta}
 \nc{\Om}{\Omega}
 \nc{\s}{\sigma}
 \nc{\ro}{\rho}
 \nc{\te}{\theta}
 \nc{\SLR}{SL_2(\br)}
 \nc{\GLR}{GL_2(\br)}
 \nc{\PGLR}{PGL_2(\br)}
 \nc{\PSLR}{PSL_2(\br)}
 \nc{\SLC}{SL(2,\bc)}
 \nc{\uH}{\mathbb{H}}
 \nc{\CH}{\mathcal{H}}
 \nc{\ck}{\mathcal{K}}
 \nc{\fD}{\mathcal{D}}
 \nc{\fE}{\mathcal{E}}
 \nc{\fO}{\mathcal{O}}
 \nc{\fX}{\mathcal{X}}
 \nc{\haf}{\frac{1}{2}}
 \nc{\qtr}{\frac{1}{4}}
 \nc{\shaf}{{\scriptstyle\frac{1}{2}}}
 \nc{\hlm}{{\scriptstyle\frac{\lambda}{2}}}
 \nc{\E}{\mathbb{E}}
 \nc{\hafc}{\frac{c}{2}}
 \nc{\8}{\infty}
 \nc{\7}{{-\infty}}
 \nc{\inv}{^{-1}}
 \nc{\eps}{\varepsilon}

\begin{document}
\currannalsline{0}{2010}

\title{Subconvexity bounds for triple $L$-functions\\ and representation theory}
\shorttitle{Subconvexity of triple $L$-functions}
\vskip 0.5cm
\rightline{\it In memory of  Ilya Piatetsky-Shapiro.}
\vskip 0.5cm

\twoauthors{Joseph Bernstein}{Andre Reznikov}
\institution{Tel Aviv University, Ramat Aviv, Israel, Bar Ilan University, Ramat-Gan, Israel}
\email{bernstei@math.tau.ac.il}
\email{reznikov@math.biu.ac.il}





\begin{abstract} We describe a new method to estimate the
trilinear period on automorphic representations of $\PGLR$. Such a
period gives rise to a special value of the triple $L$-function.
We prove a bound for the triple period which amounts to a
subconvexity bound for the corresponding special value of the
triple $L$-function. Our method is based on the study of the
analytic structure of the corresponding unique trilinear
functional on unitary representations of $\PGLR$.
\end{abstract}
\section{Introduction}
\label{intro}

\subsection{Maass forms} Let $\uH$ denote the upper half plane
equipped with the standard Riemannian metric of constant curvature
$-1$. We denote by $dv$ the associated volume element and by
$\Delta$ the corresponding Laplace-Beltrami operator on $\uH$.

Fix a discrete group $\G$ of motions of $\uH$ and consider the
Riemann surface $Y = \G \backslash \uH$. For simplicity we assume
that $Y$ is compact (the case of $Y$ of finite volume is discussed
at the end of the introduction). According to the uniformization
theorem, any compact Riemann surface $Y$ with the metric of
constant curvature $-1$ is a special case of this construction.

Consider the spectral decomposition of the operator $\Delta$ in the space
 $L^2(Y,dv)$ of functions on $Y$. It is known that the operator $\Delta$
 is non-negative and has a purely discrete spectrum;
 we will denote the
eigenvalues of $\Delta$ by $0=\mu_0< \mu_1 \leq \mu_2 \leq ...\ $.

For these eigenvalues, we always use a natural (from the
representation-theoretic point of view) parametrization
$\mu_i=\frac{1-\lm_i^2}{4}$, where $\lm_i\in \bc$. We denote by
$\phi_i=\phi_{\lm_i}$ the corresponding eigenfunctions (normalized
to have $L^2$-norm one).

In the theory of automorphic forms, the functions $\phi_{\lm_i}$
are called automorphic functions or {\it Maass forms} (after H.
Maass, \cite{M}). The study of Maass forms plays an important role
in analytic number theory, analysis and geometry. We are
interested in their analytic properties and will present a new
method of bounding some important quantities arising from
functions $\phi_{i}$.

A specific problem that we are going to address in this paper belongs
to an active area of research in the theory of automorphic
functions that studies an interplay between periods, special values of
automorphic $L$-functions and representation theory. One of the
central features of this interplay is the uniqueness of invariant
functionals associated to corresponding periods. The discovery of
this interplay goes back to classical works of E. Hecke and H.
Maass.

It is well-known that  uniqueness plays a central role in the
modern theory of automorphic functions (see \cite{PS}). The impact
that uniqueness has on the analytic behavior of periods and
$L$-functions is yet another manifestation of this principle.

\subsection{Triple products}\label{triple} \ For any three Maass
forms $\phi_i,\ \phi_j, \ \phi_k$, we define the following {\it
triple product} or {\it triple period}:\label{cijk}
\begin{eqnarray*}
c_{ijk}=\int_Y\phi_i\phi_j\phi_kdv \ .
\end{eqnarray*}

We would like to estimate the coefficient $c_{ijk}$ as a function
of parameters $\lm_i,\ \lm_j,\ \lm_k$. In particular, we would
like to find bounds for these coefficients as one or more of the
indices $i,\ j,\ k$ tend to infinity.

The bounds on the coefficient $c_{ijk}$ are related to bounds on
automorphic $L$-functions as can be seen from the following
beautiful formula of T. Watson (see \cite{Wa}, \cite{Ic}):
\begin{eqnarray}\label{wats}
\left|\int_Y\phi_i\phi_j\phi_kdv\right|^2=
\frac{\Lambda(1/2,\phi_i\otimes\phi_j\otimes\phi_k)} {
\Lambda(1,\phi_i,Ad)\Lambda(1,\phi_j,Ad)\Lambda(1,\phi_k,Ad)}\ .
\end{eqnarray}
Here the $\phi_t$ are the so-called cuspidal Hecke-Maass functions
of norm one on the Riemann surface $Y=\G\sm \uH$ arising from the
full modular group $\G=SL_2(\bz)$ or from the group of units of a
quaternion algebra. The functions
$\Lambda(s,\phi_i\otimes\phi_j\otimes\phi_k)$ and
$\Lambda(s,\phi,Ad)$ are appropriate {\it completed} automorphic
$L$-functions.

It was first discovered by R.~Rankin \cite{Ra} and A.~Selberg
\cite{Se} that the special case of above mentioned triple product  gives
rise to an automorphic $L$-function (namely, they considered the
case where one of the Maass forms is replaced by an Eisenstein
series). That allowed them to obtain analytic continuation and
effective bounds for these $L$-functions and, as an application,
to obtain one of the first non-trivial bounds for Fourier
coefficients of cusp forms towards Ramanujan's conjecture. The
relation (\ref{wats}) can be viewed as a far reaching
generalization of the original Rankin-Selberg formula. The
relation \eqref{wats} was motivated by the work  of
M.~Harris and S.~Kudla (\cite{HK}) on a conjecture of H.~Jacquet.

\subsection{Results}\label{results} \ In this paper we consider
the following problem. We fix two Maass forms $\phi = \phi_\tau$
and $ \phi'= \phi_{\tau'}$ as above and consider the coefficients
defined by the triple period:
\begin{eqnarray}\label{ci} c_i=\int_Y\phi\phi'\phi_idv\ \ \
\end{eqnarray}
as the $\phi_i$ run over an orthonormal basis of Maass forms.

Thus we see from (\ref{wats}) that the estimates of the
coefficients $c_i$ are essentially equivalent to the estimates of
the corresponding $L$-functions. One would like to have a general
method of estimating the coefficients $c_i$ and similar
quantities. This problem was raised by Selberg in his celebrated
paper \cite{Se}.

The first non-trivial observation is that the coefficients $c_i$
have exponential decay in $|\lm_i|$ as $i\to\8$. Namely, as we
have shown in \cite{BR3}, it is natural to introduce the normalized
coefficients
\begin{eqnarray}\label{d-i} d_i=\g(\lm_i)|c_i|^2\ .
\end{eqnarray}
Here $\g(\lm)$ is given by an explicit rational expression in
terms of the standard Euler $\G$-function (see \cite{BR3}) and,
for purely imaginary $\lm$, $|\lm|\to\8$, it has an asymptotic
$\g(\lm)\sim\beta |\lm|^{2}\exp(\frac{\pi}{2}|\lm|)$ with some
explicit $\beta>0$. It turns out that the normalized coefficients
$d_i$ have at most {\it polynomial growth} in $|\lm_i|$, and hence
the coefficients $c_i$ decay exponentially. This is consistent
with (\ref{wats}) and general experience from the analytic theory
of automorphic $L$-functions (see \cite{BR3}, \cite{Wa}). In
Section \ref{herm-forms} we explain a more conceptual way to
introduce the coefficients $d_i$ which is based on considerations
from representation theory.

In \cite{BR3} we proved the following mean value bound
\begin{eqnarray}\label{mean-value-convex}
\sum_{|\lm_i|\leq T}d_i\leq AT^2\ ,
\end{eqnarray}
for arbitrary $T > 1$ and some effectively computable constant
$A$.

The constant $A$ depends on the geometry of $\G$ and on parameters $\tau$,
$\tau'$ of eigenfunctions $\phi$, $\phi'$.

According to Weyl's law for the spectrum of the
Laplace-Beltrami operator $\Dl$ on $Y$, the number of terms in this
sum is of order $C T^2$. So this formula says that on average the
coefficients $d_i$ are bounded by some constant.

More precisely, let us fix an interval $I \subset \br$ centered at the point $T$
and consider the finite set of all Maass forms $\phi_i$ with
parameter $|\lm_i|$ inside this interval. Then the average value
of coefficients $d_i$ in this set is bounded by a constant {\it
provided} the interval $I$ is long enough (i.e., of size $\approx
T$).

Note that the best individual bound which we can get from this
formula is $d_i\leq A|\lm_i|^2$. For Hecke-Maass forms this bound
corresponds to the convexity bound for the corresponding
$L$-function via Watson formula (\ref{wats}).

The central result of this paper is the bound for the sum of the coefficients
$d_i$ over a shorter interval. Namely, we prove the following

\begin{thm}{thm}There exist effectively computable constants $B,\ b >0$ such that,
for an arbitrary $T >1$, we have the following bound
\begin{eqnarray}\label{sub-convex-thm}
\sum_{|\lm_i| \in I_T} d_i\leq B T^{\frac{5}{3}}\ ,
\end{eqnarray}
where $I_T$ is the interval of size $b T^{1/3}$ centered at $T$.
\end{thm}

The exponent $5/3$ above appears for the reason similar to the
appearance of the exponent $1/3$ in the asymptotic of the Airy
integral (namely, a degenerate critical point in the phase of an
oscillatory integral; see Remark \ref{rem1}).

The constant $B$ depends on the geometry of $X$ and on parameters
$\tau$, $\tau'$ (see Remark~\ref{D-depend}). The constant $b$ depends on parameters $\tau$,
$\tau'$ only.

Note that the theorem gives an individual bound $d_i \leq B
|\lm_i|^{\frac{5}{3}}$ (for $|\lm_i| > 1$). Thanks to the Watson formula
(\ref{wats}) and a lower bound of H. Iwaniec
$|L(1,\phi_{\lm_i},Ad)|\gg |\lm_i|^{-\eps}$ (see \cite{I}), this
leads to the following {\it {subconvexity}} bound for the triple
$L$-function (for more on the relation between triple period and
special values of $L$-functions, see \cite{Wa}, \cite{Ic}).

\begin{cor}{cor}Let $\phi$ and $\phi'$ be fixed Hecke-Maass cusp forms. For any
$\eps>0$,
 there exists $C_\eps>0$
such that the bound
\begin{eqnarray}\label{sub} L(1/2\ ,
\phi\otimes\phi'\otimes\phi_{\lm_i})\leq C_\eps|\lm_i|^{{\frac{5}{3}}+\eps}
\end{eqnarray} holds for any
Hecke-Maass form $\phi_{\lm_i}$.
\end{cor}

The convexity bound for the triple $L$-function corresponds to
\eqref{sub} with the exponent $5/3$ replaced by $2$. We refer to
\cite{IS} for a discussion of the subconvexity problem which is at
the core of modern analytic number theory. We note that the above
bound is the first subconvexity bound for an $L$-function of
degree $8$ which does not split in a product of smaller degree $L$-functions.
 All previous subconvexity results were obtained for
$L$-functions of degree at most $4$.

In  \cite{V} A. Venkatesh obtained a subconvexity
bound for the triple $L$-function in the level aspect (i.e., with
respect to a tower of congruence subgroups $\G(N)$ as $N\to\8$).
His method is quite different from the method we present in this
paper and is based on ideas from ergodic theory.

We formulate a natural

\begin{conj}{conj}For any $\eps>0$ we have $d_i\ll|\lm_i|^{\eps}$ .
\end{conj}

For Hecke-Maass forms on congruence subgroups, this conjecture is
consistent with the Lindel\"{o}f conjecture for the triple
$L$-functions (for more details, see \cite{BR3} and \cite{Wa}).

\subsubsection{Remarks.}\ 1. Our results can be extended to the
case of a general finite co-volume lattice $\G \subset G$ (see Remark \ref{proof-remks} for more detail).

2. First results on the exact exponential decay of triple products
for a general lattice $\G$ and holomorphic forms were obtained by A.
Good \cite {Go} using Poincar\'{e} series.
P. Sarnak \cite{Sa1} discovered ingenious analytic continuation of
Maass forms to the complexification of the Riemann surface $Y$  to obtain somewhat
weaker results for Maass forms (for representation-theoretic approach to this method and
generalizations, see \cite{BR1} and \cite{KS}). Our present
method seems to be completely different and avoids analytic
continuation.

3. We would like to stress that the bound for the triple product
in Theorem \ref{thm} is valid for a general lattice $\G$,
including {\it non-arithmetic} lattices. In fact, in our method we
do not use Fourier coefficients or Hecke eigenvalues through which
one usually accesses values of $L$-functions for congruence
subgroups. Our method gives estimates for {\it periods} of
automorphic functions directly and $L$-functions appear only
through the Watson formula \eqref{wats} (the same is true for the
method of Venkatesh \cite{V}).
\\
\\

The paper is organized as follows. The next section is devoted
to a detailed explanation of ideas behind the method of the proof
of Theorem \ref{thm}. The main body of the paper (Sections 3-10) is
devoted to the proof. Two Appendices containing technical calculations
conclude the paper.
The numbering in the paper is organized as follows. Each subsection has a unique
Theorem, Proposition, Lemma etc., and these are numbered by the corresponding section.
Equations are numbered continuously within each section.

\subsection*{Acknowledgements}
It is a pleasure to thank Peter Sarnak for stimulating discussions
and support of this work. We would like to thank Herv\'{e} Jacquet
for a valuable comment, Misha Sodin for  analytic advice, and the
referee for constructive comments.

The research was partially supported by a BSF grant, by a GIF grant, by
the Excellency Center ``Application of Algebraic Geometry and
Logic to Representation theory'' of the Israel Science
Foundation and by the Emmy Noether Institute for Mathematics (the Center
of Minerva Foundation of Germany). The
results of this paper were mostly obtained during our visits to the
Max-Planck Institute in Bonn and in Leipzig, to the Courant
Institute, and to Weizmann Institute. We would like to thank these
institutions for their excellent atmosphere.

\section {Outline of the proof}
We describe now the general ideas behind our proof. It is based on
ideas from representation theory (for a detailed account of the
corresponding setting, see \cite{BR3} and Section~\ref{3-prod}
below). In what follows we sketch the method of the proof whose
technical details appear in the rest of the paper.

\subsection{Automorphic representations}{\label{automorphic}
 \ Let $G$ denote the group of all motions of $\uH$. This group is naturally
isomorphic to $\PGLR$ and as a $G$-space $\uH$ is naturally
isomorphic to $G / K$, where $K = PO(2)$ is the standard maximal
compact subgroup of $G$.

By definition, $\G$ is a subgroup of $G$. The space $X = \G
\backslash G$ with the natural right action of $G$ is called an
{\it automorphic space}. We will identify the Riemann surface $Y =
\G \backslash \uH$ with $X / K=\G\setminus G/K$.

We use the standard language of automorphic representations (see
\cite{GGPS} and Section \ref{reps} below). Let $(\pi,G,V)$ be an
irreducible smooth representation of $G$. An {\it automorphic
structure} on $V$ is a continuous $G$-morphism $\nu:V\to C^\8(X)$.

The pair $(\pi,\nu)$ consisting of an abstract representation
$(\pi,V)$ and the automorphic structure $\nu$ will be called an
{\it automorphic representation}. This terminology is slightly
more precise then the standard one. We find it more convenient for
our purposes.

We always assume that $(\pi, V)$
is unitary (i.e., $V$ is equipped with a positive definite $G$-invariant
Hermitian form $P$), and that the automorphic structure $\nu$ is
compatible with the invariant Hermitian form $P$.

We will usually present the abstract representation $(\pi, V)$ by an
explicit model. We will deal mostly with class one irreducible
representations of $G$ (i.e., those with a non-zero $K$-fixed
vector). If $(\pi,V)$ is a non-trivial class one representation we
use for it the model $V = V_\lm$, where $\lm\in i\br\cup(0,1)$ and
$V_\lm$ is the space of smooth even homogeneous functions on
$\br^2\sm 0$ of the homogeneous degree $\lm-1$ (see \cite{G5},
\cite{BR3}). We denote by $e_\lm\in V_\lm$ the function taking
constant value $1$ on $S^1\subset \br^2\sm 0$. This gives a
$K$-invariant vector in the representation $V_\lm$ which we call the standard $K$-fixed vector in $V_\lm$.
We normalize the invariant Hermitian form $P$ on $V_\lm$ by the condition $P(e_\lm)=1$.

The theorem of Gelfand and Fomin states that all Maass forms (or more
generally automorphic functions) could be obtained as special
vectors in appropriate automorphic representations (see
\cite{GGPS}). Namely,  a Maass form $\phi=\nu(e_\lm)$ corresponding to an
automorphic structure $\nu$ on a representation with a model $V_\lm$ has the eigenvalue
$\mu=\frac{1-\lm^2}{4}$.

We translate various questions about Maass forms into
corresponding questions about associated automorphic
representations. This allows us to employ powerful methods of
representation theory.

\subsection{}\label{uniq-introd} Let us fix two (nontrivial) automorphic
representations $(\pi,\nu)$ and $(\pi',\nu')$. We assume that both
are representations of class one (i.e., $V\simeq V_\tau$ and
$V'\simeq V_{\tau'}$, $\tau,\ \tau'\in i\br\cup(0,1)$). These give
rise to Maass forms $\phi=\nu(e_\tau)$ and
$\phi'=\nu '(e_{\tau'})$. Let $(\pi_i,V_i, \nu_i)$ be a third
automorphic representation (which we a going to vary) with the
parameter $\lm_i$ (i.e., $V_i\simeq V_{\lm_i}$).

The triple product $c_i=\int_Y\phi\phi'\phi_idv$ extends to a
$G$-equivariant trilinear functional on the corresponding
automorphic representations $l^{aut}_i:V\otimes V'\otimes
V_i\to\bc$.

Next we use a general result from representation theory that such a
$G$-equivariant trilinear functional is unique up to a scalar,
i.e., that $\dim\Mor_G(V\otimes V'\otimes V'',\bc)\leq 1$ for any
smooth irreducible representations $V,\ V',\ V''$ of $G$ (see
\cite{O}, \cite{Pr}, \cite{Lo} and the discussion in \cite{BR3}).
This implies that the automorphic functional $l^{aut}_i$ is
proportional to some explicit {\it model} functional
$l^{mod}_{\lm_i}$. In \cite{BR3} we gave a description of such
a model functional $l^{mod}_{\lm}: V \otimes V' \otimes V_{\lm} \to
\bc$ for any $\lm$ using explicit realizations of representations
$V$, $V'$ and $V_{\lm}$ of the group $G$ in spaces of homogeneous
functions; it is important that the model functional knows nothing
about the automorphic picture and carries no arithmetic information.

Thus we can write $l^{aut}_i = a_i \cdot l^{mod}_{\lambda_i}$ for
some constant $a_i$, and hence
\begin{eqnarray}\label{c-a} c_i=
l^{aut}_i(e_\tau\otimes e_{\tau'} \otimes e_{\lambda_i}) = a_i
\cdot l^{mod}_{\lm_i}(e_\tau \otimes e_{\tau'} \otimes
e_{\lambda_i})\ ,
\end{eqnarray}
where $e_\tau,\ e_{\tau'},\ e_{\lambda_i}$ are the standard  K-invariant unit
vectors in representations $V, V'$ and $V_{\lambda_i}$
corresponding to the automorphic forms $\phi$, $\phi'$ and
$\phi_i$.

It turns out that the proportionality coefficient $ a_i$ in
\eqref{c-a} carries  important \lq\lq automorphic" information
while the second factor carries no arithmetic information and can
be computed in terms of $\G$-functions using explicit realizations
of representations $V_\tau$, $V_{\tau'}$ and $V_{\lambda}$ (see
Appendix in \cite{BR3} where this computation is carried out).
This second factor is responsible for the exponential decay, while
the first factor $a_i$ has a polynomial behavior in parameter
$\lm_i$. An explicit computation shows (see loc. cit.) that
$|c_i|^2=\frac{1}{\g(\lm_i)}|a_i|^2$, and hence $d_i=|a_i|^2$
(where the function $\g(\lm)$ was described in Section
\ref{results}).

So, from now on we will deal with
coefficients $d_i$ and no longer refer to coefficients $a_i$ and $c_i$ at all.

\subsection{Hermitian forms}\label{hermitian} \ In order to
estimate the quantities $d_i$, we consider the space $E = V_\tau
\otimes V_{\tau'}$ and use the fact that the coefficients $d_i$
appear in the spectral decomposition of the following {\it
geometrically defined} non-negative Hermitian form $H_\Dl$ on $E$
(for a detailed discussion, see \cite{BR3}).

Consider the space $C^\8(X\times X)$. The diagonal $\Dl:X \to
X\times X$ gives rise to the restriction morphism $r_\Dl :
C^\8(X\times X)\to C^\8(X)$. We define a non-negative Hermitian
form $H_\Dl$ on $C^\8(X\times X)$ by setting $H_\Dl =
(r_\Dl)^{\ast}(P_X)$, where $P_X$ is the standard $L^2$ Hermitian
form on $C^\8(X)$ i.e.,
$$H_\Dl(w)= P_X(r_\Dl(w)) = \int_{X}|r_\Dl(w)|^2d\mu_X $$
 for any $w\in C^\8(X\times X)$.
We call the restriction of the Hermitian form $H_\Dl$ to the
subspace $E\subset C^\8(X\times X)$ the {\it diagonal} Hermitian
form and denote it by the same letter.

We will describe the spectral decomposition of the Hermitian form
$H_\Dl$ in terms of Hermitian forms corresponding to trilinear
functionals. Namely, if $L$ is a pre-unitary representation of $G$
with $G$-invariant Hermitian norm $|| \ ||_L$, then every $G$-invariant
trilinear functional $l: V \otimes V' \otimes L \to \bc$ defines a Hermitian form $H^l$ on $E$ by $H^l(w)
=\sup\limits_{||u||_L = 1} |l(w \otimes u)|^2 \ $.

Here is another description of this form (see
\cite{BR3}). The functional $l: V \otimes V' \otimes L \to \bc$ gives
rise to a $G$-intertwining morphism $T^l: E \to L^*$ which image
lies in the smooth part $\tilde L^*$ of $L^*$. Then the form $H^l$ is just the
pull back of the Hermitian form on $\tilde L^*$ corresponding to the
inner product on $L$.

Consider the orthogonal decomposition $L^2(X) =\left(\bigoplus_i
V_i\right)\oplus\left( \bigoplus_{\kappa} V_{\kappa}\right)$ where
$V_i$ correspond to Maass forms and $V_{\kappa}$ correspond to
representations of discrete series. Every $G$-invariant subspace
$L \subset L^2(X)$ defines a trilinear functional $l :E \otimes L
\to \bc$ and hence a Hermitian form $H^l$ on $E$. Hence, the
decomposition of $L^2(X)$ gives rise to the corresponding
decomposition $$H_\Dl = \sum H_i^{aut} + \sum H_{\kappa}^{aut}$$
of Hermitian forms (see \cite{BR3}).

We denote by $H_\lm$ the {\it model} Hermitian form corresponding
to the {\it model} trilinear functional $l^{mod}_\lm :V \otimes V'
\otimes V_\lm \to \bc$. The uniqueness of trilinear functionals
mentioned in Section \ref{uniq-introd} (i.e., the formula \eqref{c-a}) implies that $H_i^{aut} =
d_i H_{\lm_i}$. This leads us to

{\bf The basic spectral identity}
\begin{eqnarray}\label{besseli} H_\Dl\ = \sum_i d_iH_{\lm_i} + \sum_{\kap}
H_{\kap}^{aut},
\end{eqnarray}

Of course, one can introduce similar model trilinear functionals
for the discrete series representations $V_\kap$ and the
corresponding coefficients $d_\kap$ via $H_{\kap}^{aut}=d_\kap
H_\kap$. We will not need these in this paper (in fact, in this paper we are
trying to avoid computations with the discrete series
representations; see Remark \ref{rem-trip}).

We will mostly use the fact that for every vector $w \in E$ this
basic spectral identity gives us an inequality
\begin{eqnarray}\label{bessel} \sum_i d_i H_{\lm_i}(w) \leq H_\Dl (w)
\end{eqnarray}
which turns into an {\it equality} if the vector $r_\Dl(w)$ does
not have projection to discrete series representations (for
example, if the vector $w$ is {\it invariant} with respect to the
diagonal action of $K$ on $E$).

We can use this inequality to bound coefficients $d_i$. Namely,
for a given vector $w \in E$ we usually can compute the values of
the weight function $H_{\lm}(w)$ by explicit computations in the
model of representations $V,\ V',\ V_{\lambda}$. It is usually
much more difficult to get reasonable estimates of the right hand
side $H_{\Dl}(w)$ since it refers to the automorphic picture. In
cases when we manage to do this we get some bounds for the
coefficients $d_i$.

\subsection{Mean-value estimates} \ In \cite{BR3}, using the
geometric properties of the diagonal form and explicit estimates
of forms $H_\lm$, we established the mean-value bound
(\ref{mean-value-convex}):
$$\sum\limits_{|\lm_i| \leq T} d_i \leq A T^2\ .$$
Roughly speaking, the proof of this bound is based on
the fact that while the value of the form $H_{\Dl}$ on a given
vector $w \in E$ is very difficult to control, we can show that
for many vectors $w$ the value $H_\Dl(w)$ can be bounded by
$P_E(w)$, where $P_E$ is the Hermitian form which defines the
standard unitary structure on $E$.

More precisely, consider the natural representation $\sigma = \pi
\otimes \pi'$ of the group $G\times G$ on the space $E$. Then for
a given compact neighborhood $U \subset G \times G$ of the
identity element, there exists a constant $C$ such that for any
vector $w \in E$, the inequality $H_{\Dl}(\sigma(g)w) \leq C
P_E(w)$ holds for at least half of the points $g \in U$. This
follows from the fact that the average over $U$ of the quantity
$H_{\Dl}(\sigma(g)w)$ is bounded by $C P_E(w)/2$.

This allows us, for every $T\geq 1$, to show the {\it existence} of  a vector $w \in E$ such
that $H_{\Dl}(w) \leq C T^2$ and $H_\lm(w) \geq c$ for all $\lm$ satisfying $|\lm| \leq T$. The bound \eqref{bessel} then implies the mean-value bound \eqref{mean-value-convex}.

\subsection{Bounds for sums over shorter intervals} \
 The main starting point of our approach to the
subconvexity bound is the inequality (\ref{bessel}) for Hermitian
forms. For a given $T>1$, we construct a test vector $w_T\in E$
such that the weight function $\lm\mapsto H_\lm(w_T)$ has a sharp
peak near $|\lm|=T$ (i.e., a vector satisfying the condition
\eqref{lower-bd} below).

The problem is how to estimate effectively $H_\Dl(w_T)$. The idea
is that the Hermitian form $H_\Dl$ is geometrically defined and,
as a result, satisfies some non-trivial bounds, symmetries, etc.
None of the explicit {\it model} Hermitian forms $H_\lm$ satisfies
similar properties. By applying these symmetries to the vector
$w_T$, we construct a new vector $\tilde w_T$ and from the geometry
of the automorphic space $X$, we deduce the bound $H_\Dl(w_T)\leq
H_\Dl(\tilde{w}_T)$.

On the other hand, the weight function $H_\lm(\tilde w_T)$ in the
spectral decomposition $H_\Dl(\tilde w_T) = \sum d_i H_{\lm_i}(
\tilde w_T)$ for $\tilde w_T$ behaves quite differently from the
weight function $H_\lm(w_T)$ for $w_T$. Namely, the function
$H_\lm(\tilde{w}_T)$ behaves regularly (i.e., satisfies condition
(\ref{(ii)}) below), while the weight function $H_\lm(w_T)$ has a
sharp peak near $|\lm| = T$.

The regularity of the function $H_\lm(\tilde w_T)$ coupled with
the mean-value bound (\ref{mean-value-convex}) allows us to prove
a sharp upper bound on the value of $H_\Dl(\tilde{w}_T)$ by purely
spectral considerations (in the cases that we consider there is no
contribution from discrete series). We do not see how to get such
sharp bound by geometric considerations working on the automorphic
space $X\times X$.

Using this bound for $H_\Dl(\tilde w_T)$ and the inequality
$H_\Dl(w_T)\leq H_\Dl(\tilde{w}_T)$, we obtain a non-trivial bound
for $H_\Dl(w_T)$ and, as a result, the desired bound for the
coefficients $d_i$.

We now describe this strategy in more detail.

\subsection{Proof of Theorem \ref{thm}}\label{test-vect-sect}
We only consider the case of representations of the principal
series, i.e., we assume that $V = V_{\tau}$, $V' =
V_{\tau'}$ for some $ \tau,\ \tau' \in i \mathbb{R}$; \ the case of
representations of the complementary series can be treated
similarly.

We denote by $\nu$ and $\nu'$ the corresponding automorphic
realizations of $V$ and $V'$. We choose an orthonormal basis
$\{e_n\}_{n\in 2\bz}$ in $V$ consisting of $K$-types and similarly
an orthonormal basis $\{e'_n\}$ in $V'$.

Vectors $w_n =e_n\otimes
 e'_{-n} \in E = V \otimes V'$ will play an important role in our computations.

Let us set
\begin{eqnarray}\label{def-of-S} \mathcal{S}=2(|\tau|+|\tau'|)+1
\end{eqnarray} the constant depending  on parameters of representations
$V$ and $V'$ only. For a given $T\geq \mathcal{S}$, we choose an even integer $n$ such that
$|T-2n|\leq 10$ and set
\begin{eqnarray}\label{test-vectors1} w_T=w_n=e_n\otimes e'_{-n}\ .
\end{eqnarray}
In fact, all we need is that $|T-2n|$ remain bounded as $T\to\8$.

By a direct computation involving stationary phase method, we show in
Section \ref{proof-i} that the following lower bound holds

{\bf First spectral bound:}

There exist constants $b, c > 0$ such that
\begin{eqnarray}\label{lower-bd} \
 H_\lm(w_T)\geq c \ T^{-5/3}\text{ for } |\lm| \in
I_T\ ,
\end{eqnarray}
where $I_T$ is the interval of length $bT^{1/3}$ centered at the point
$T$.

This inequality together with the bound $\sum_i d_i H_{\lm_i}(w_T)
\leq H_\Dl (w_T)$ (see \eqref{bessel}) imply the bound
\begin{eqnarray}\label{upper-bd2}
 \sum\limits_{|\lm_i|\in I_T}d_i\leq C T^{5/3}H_\Dl(w_T)\ ,
\end{eqnarray} for some constant $C$.

Now we claim that the quantity $H_\Dl(w_T)$ is uniformly bounded by
some constant $D$ which does not depend on $T$.
 Namely we can write \begin{eqnarray*}
\ H_\Dl(w_T)=\int\limits_{X}|\nu(e_n)|^2|\nu'(e'_{-n})|^2\
 d\mu_X\leq  \haf\left(||\nu(e_n)||^4_{L^4(X)}+||\nu'(e'_{-n})||^4_{L^4(X)}\right)\ .\nonumber
\end{eqnarray*}

Hence the necessary bound follows from the following result which,
we feel, is of
independent interest.

\begin{thm}{thm2}For a fixed class
one automorphic representation $\nu:V\to C^\8(X)$, there exists a
constant $D>0$ such that $||\nu(e_n)||_{L^4(X)}\leq D$ for all $n$.
\end{thm}

This finishes the proof of Theorem \ref{thm}.\qed

\begin{rem}{}One would expect that $L^4$-norms of $K$-types for
representations of
the discrete series are uniformly bounded as well. It is a very
interesting and deep question to
study dependence of the constant $D$ in Theorem \ref{thm2}
on the parameter $\tau$ of the
automorphic representation and on the subgroup $\G$ (for a discussion, see Remark \ref{D-depend}). Moreover, it would be interesting to identify (as a norm on an abstract representation $\pi_\tau$) the $G$-{\it invariant} (non-Hermitian) norm which  the $L^4$-norm on $X$ induces on the representation $\pi_\tau$ via automorphic isometry $\nu_\tau$ (see a discussion in \cite{BR1}).

Another interesting question is an analog of the above theorem for a cuspidal representation for a non-uniform $\Gamma$. Specifically, we would like to know whether  $L^4$-norm of $K$-types are uniformly bounded for a fixed cuspidal representation (compare to Remark 2, Section \ref{proof-remks}).
\end{rem}

\subsection{$L^4$-norms of $K$-types} We now explain the proof of
the uniform bound for $L^4$-norm of $K$-types (i.e., Theorem \ref{thm2}).

Let $\bar V$ be the complex conjugate to $V$ representation. The representation $\bar V$ is also an automorphic
representation with the realization $\bar\nu:\bar V\to C^\8(X)$
(see details in Section \ref{complex-c}). For the proof of Theorem \ref{thm2} it is enough to consider the
setup described above (i.e., the space $E$, forms $H_\Dl$, $H_\lm$,
etc.) for the special case when $V'$ is isomorphic to the representation $\bar V$.

We only consider the case of representations of the principal
series, i.e., we assume that $V = V_{\tau}$ and $V' =\bar V=
V_{-\tau}$ for some $ \tau \in i \mathbb{R}$; \ the case of
representations of the complementary series can be treated
similarly.

Choose an orthonormal basis $\{e_n\}_{n\in 2\bz}$ in $V$
consisting of $K$-types.
We denote by $\{e'_n=\overline{ e_{-n}}=c(e_{-n})\}$ the complex
conjugate basis in $\bar V$ (note that $e'_n$ is of the $K$-type
$n$).

For a given $n\in 2\bz$, we set
\begin{eqnarray}\label{test-vectors} w_n=e_n\otimes e'_{-n}\ \ \ \text{
and}\ \ \ \ \tilde w_n= w_n + w_{n+2} .
\end{eqnarray}
With such a choice of test vectors we have the following bounds.

{\bf Geometric bound:}
\begin{eqnarray}\label{geom-bound}
\label{(star)} H_\Dl(w_n)\leq H_\Dl(\tilde w_n)\ .
\end{eqnarray}

{\bf Second spectral bound:}

There exists a constant $C'$ such that
\begin{eqnarray}\label{(ii)}\hskip 1cm H_\lm(\tilde w_n)\leq
\begin{cases}
C'(1+|n|)\inv|\lm|^{-1}+C'|\lm|^{-3}&\text{for all}\ \mathcal{S}\leq|\lm|\leq 4|n|\ ,\\
 C'|\lm|^{-3}&\text{for all}\ |\lm|>4|n|\ . \\
\end{cases}
\end{eqnarray}
Here $\mathcal{S}$ is as in \eqref{def-of-S}.

Using the bound \eqref{(ii)} we will get the following {\it sharp}
estimate of $H_\Dl(\tilde w_n)$ (see Proposition \ref{main-prop}):
\begin{eqnarray}\label{w-leq-D}
H_\Dl(\tilde w_n)\leq D\ \
\end{eqnarray}
with some explicit constant $D>0$ (for the proof, see Section
\ref{proof-geom-bd}). Bounds \eqref{w-leq-D} and \eqref{geom-bound}
imply the bound for the $L^4$-norm of $K$-types since in this case
 $H_\Dl(w_n)=||\nu(e_n)||_{L^4}^4$.

The bound \eqref{w-leq-D} follows from the identity
$H_\Dl(\tilde w) = \sum d_i H_{\lm_i}(\tilde w)$ (see
\eqref{besseli}), the spectral bound \eqref{(ii)} and the
mean-value bound \eqref{mean-value-convex} for the coefficients
$d_i$. The low spectrum contribution for $|\lm_i|\leq \mathcal{S}$  is bounded by an argument based on the Sobolev restriction theorem (see Section \ref{low-spec}) . We also use the fact that there are no contribution to
$H_\Dl(\tilde w)$ coming from the discrete series since the vector
$\tilde w$ is $\Delta K$-invariant.

\subsubsection{Proof of the geometric bound (\ref{geom-bound})}
\label{g-bound-proof} \ The inequality (\ref{(star)}) easily follows
from the pointwise bound on $X$. Namely, in the automorphic
realization, the vector $w_n = e_n\otimes e'_{-n}$ is represented by
a function whose restriction $u_n = r_{\Delta}(\nu_E(w_n))$ to the
diagonal is {\it non-negative} (see also Section \ref{aut-positive})
$$
u_n(x) = \nu(e_n)(x)\cdot\bar\nu(e'_{-n})(x)=|\nu(e_n)(x)|^2\geq
0.
$$
From this we see that
\begin{eqnarray*}\label{w-w-bound}
H_\Dl(w_T)=\int_{X}|u_n(x)|^2d\mu_X \leq
\int_{X}|u_n(x)+u_{n+2}(x)|^2d\mu_X= H_\Dl(\tilde{w}_T)\ .\ \ \ \qed
\end{eqnarray*}

\subsubsection{Sketch of proof of the spectral bounds \eqref{lower-bd}
and \eqref{(ii)}.}\label{sketch} Proof of these bounds is carried out
by the standard application of the stationary phase method and the
Van der Corput lemma. It constitutes the main {\it technical} bulk
of the paper. We will use the explicit form of the kernel defining
Hermitian forms $H_\lm$ in the model realizations of
representations $V$, $V'$ and $V_\lm$. Namely, we use the standard
realization of these representations in the space
$C^\8_{even}(S^1)$ of even functions on $S^1$ (see \cite{BR3} and
Section \ref{automorphic}). Under this identification, the basis
$\{e_n\}_{n\in 2\bz}$ becomes the standard basis of exponents
$\{e_n=e^{int}\}$, where $0\leq t<2\pi$ is the standard
parameter on $S^1$.

In \cite{BR3}, Section 5, we described how to write down an
invariant functional for principal series representations. Namely,
let $V =V_\tau$,
 $ V' = V_{\tau'}$ with $\tau,\ \tau'\in i\br$. In the circle model
of representation $V_\tau$, $V_{\tau'}$, $V_\lm$, the following
kernel on the space $V_\tau\otimes V_{\tau'}
\otimes V_\lm\simeq C^\8((S^1)^3)$ defines an
invariant functional
kernel on $(S^1)^3$:
\begin{eqnarray*}K_\lm(x,y,z)=
|\sin(x-y)|^{\frac{-1-\tau-\tau'+\lm}{2}}|\sin
(x-z)|^{\frac{-1-\tau+\tau'-\lm}{2}}|\sin
(y-z)|^{\frac{-1+\tau-\tau'-\lm}{2}}\ ,\end{eqnarray*}
where $x,\ y,\ z\in S^1$.
We denote this functional by $l^{mod}_\lm$. Using the kernel
$K_\lm(x,y,z)$, we can define the Hermitian forms $H_\lm$ on $E\simeq
C^\8(S^1\times S^1) $ by the corresponding oscillatory integral
(over $(S^1)^4$; see Section \ref {K-reduction}). This allows us to
use the stationary phase method in the proof of bounds
\eqref{lower-bd} and (\ref{(ii)}).

Here appears the main difference between test vectors $w_n$ and
$\tilde w_n$ . It manifests itself in the form of the oscillating
integrals computing $H_\lm(w_n)$ and $H_\lm(\tilde w_n)$. Namely,
both of these integrals have the same phase function which has a
{\it degenerate} critical point. The main difference between them is
that for the vector $w_T$ the corresponding integral has a non-zero
amplitude at this critical point (this gives the crucial lower bound
\eqref{lower-bd}) and for $\tilde w_T$ the amplitude vanishes at the critical point
(resulting in bounds \eqref{(ii)}).

In fact, we will use the values of $H_\lm(w)$ only for $\Dl
K$-invariant vectors $w \in E$. This considerably simplifies our
computations since we can reduce them to two repeated integrations
in one variable and use the stationary phase method in {\it one
variable}.

\begin{rems}{rem1}1. The existence of vectors satisfying spectral conditions
\eqref{lower-bd} and \eqref{(ii)} allows us to shorten the
summation over the spectrum, comparatively to the range of the
summation in the convexity bound \eqref{mean-value-convex}. This
is necessary if one wants to deduce a subconvexity bound from the
Bessel inequality of Hermitian forms \eqref{bessel} since the
convexity bound \eqref{mean-value-convex} is essentially {\it
sharp} (see \cite{Re1}). This approach to the subconvexity is
reminiscent of the classical amplification method introduced by
Selberg (see \cite{Mi}, \cite{MiV} for the review of the state of
the art subconvexity results). Usually one uses a variant of a
trace formula to control the so-called off-diagonal terms arising
after shortening the sum. In our approach there is no use of the
Selberg or the Kuznetsov trace formulas. Instead, we use the
hidden symmetries of the diagonal form $H_\Dl$.

2. The origin of our exponent $5/3=2(1/2+1/3)$ in the main Theorem
\ref{thm} (i.e., the bound \eqref{sub-convex-thm}) is directly
related to the exponent $1/3$ in the well-known properties of the
Airy function. In fact, we reduce the proof of the crucial lower
bound \eqref{lower-bd} to the asymptotic of the Airy integral (see
Proposition \ref{H-to-Airy}).

3. After obtaining results presented in this paper, we realized that
there exists another possible approach to
bounds for triple and other periods of automorphic functions.
It is based on the notion of strong
Gelfand pairs (see \cite{Gr} and references therein).
This approach is presented in \cite{Re2}.

There is one technical complication in the approach based on Gelfand pairs, though.
We where not able to
produce the desired family of test vectors which is
also $\Dl K\times \Dl K$-invariant. Without this property one has
to consider terms in the spectral decomposition
\eqref{besseli} coming from the discrete series
representations. It is more cumbersome to study model trilinear
functionals on discrete series as these representations do not
have nice geometric models. As a result, in this paper we use
another property of the form $H_\Dl$ , the extra positivity
provided by the Cauchy-Schwartz inequality (see Section
\ref{g-bound-proof}), instead of the associated Gelfand pairs
structure. We hope to return to this subject elsewhere.
\end{rems}

\section{Representation-theoretic setting}
\label{reps} \subsection{} We recall the standard connection
between Maass forms and representation theory of $\PGLR$ (see
\cite{GGPS}). Most of the material in the next three sections is
taken from \cite{BR3}, where it is discussed in more detail.

\subsubsection{Automorphic space} \label{autspace}

Let $\uH$ be the upper half plane with the hyperbolic metric of
constant curvature $-1$. The group $\SLR$ acts on $\uH$ by
fractional linear transformations. This action allows to identify
the group $\PSLR$ with the group of all orientation preserving
motions of $\uH$. For reasons explained below (see Remark
\ref{rem-uni}), we would like to work with the group $G$ of all
motions of $\uH$; this group is isomorphic to $\PGLR$. Hence
throughout the paper we consider the group $G=\PGLR$ and denote by
$K$ its standard maximal compact subgroup $K=PO(2)$. We have natural
identification $G/K=\uH$.

We fix a discrete co-compact subgroup $\G \subset G$ and set $X=\G
\sm G$. We fix the unique $G$-invariant
measure $\mu_X$ on $X$ of {\it total mass one}. The group $G$ acts on $X$ (from the right) and hence on
the space of functions on $X$.  Let
$L^2(X)=L^2(X,d\mu_X)$ be the space of square integrable functions
and $(\Pi_X, G, L^2(X))$ the corresponding unitary representation.
We will denote by $P_X$ the Hermitian form on $L^2(X)$ given by
the inner product.

\subsubsection{Automorphic representations} \label{ureps}
Let $(\pi,G,V)$ be an irreducible smooth Fr\'{e}chet
representation of $G$ (see \cite{Cass} where they are called
smooth representations of moderate growth).

\defn{aut-def}{}An {\it automorphic structure} on $(\pi,V)$
is a continuous $G$-morphism $\nu:V\to C^\8(X)$.

We call an automorphic representation a pair $(\pi,\nu)$ of a
representation and the automorphic structure on it. In this paper
we always assume that $(\pi,V)$ is irreducible, admissible and
also assume that $(\pi,V)$ is unitary. This means that $V$ is
equipped with a $G$-invariant positive definite Hermitian form $P$,
and $V$ is the space of smooth vectors in the completion $L$ of
$V$ with respect to $P$. An automorphic structure $\nu:V\to
C^\8(X)$ is assumed to be normalized, i.e., we assume that $P=\nu^*(P_X)$.

\subsubsection{Automorphic representations and Maass
forms}\label{rep-models} Let $(\pi_\lm,G,V_\lm)$ be a
representation of the generalized principal series corresponding
to $\lm\in\bc$. The space $V_\lm$ is the space of smooth even
homogeneous functions on $\br^2\sm 0$ of the homogeneous degree
$\lm-1$ (which means that $f(ax,ay)=|a|^{\lm-1}f(x,y)$ for all
$a\in\br \setminus 0$) with the action of $GL(2,\br)$ given by
$\pi_\lm (g)f(x,y)= f(g\inv (x,y))|\det g|^{(\lm-1)/2}$ (see
\cite{G5}).

In explicit computations it is often convenient to pass
from the plane model to a circle model.  Namely, the restriction of
functions in $V_\lm$ to the unit circle $S^1 \subset \br^2$ defines
an isomorphism of the space $V_\lm$ with the space
$C_{even}^\8(S^1)$ of even smooth functions on $S^1$, so we can
think about vectors in $V_\lm$ as functions on $S^1$.
The constant function $1$ on $S^1$ corresponds to the standard unit $K$-invariant vector $e_\lm\in V_\lm$.
We normalize the invariant Hermitian form $P$ by the condition $P(e_\lm)=1$. For
$\lm\in i\br$, this corresponds to the standard Hermitian form
$\langle f,g\rangle_{V_\lm}=1/2\pi\int_{S^1}f\bar gd\theta$ on (even) functions on
$S^1$.

Suppose $\nu:V\to C^\8(X)$ is an automorphic structure on $V_\lm$.
Then $\phi_\lm=\nu(e_\lm)\in C^\8(X)^K=C^\8(Y)$ is a Maass form
with the eigenvalue $\mu=\frac{1-\lm^2}{4}$.

This construction, which is due to Gelfand and Fomin, gives
a one-to-one correspondence between Maass forms and class one
automorphic representations (and more generally between
automorphic forms and automorphic representations of $G$). We refer to \cite{GGPS} for a more detailed discussion (see also \cite{BR3}).

\subsubsection{Decomposition of the representation $(\Pi_X, G,
L^2(X))$}\label{spec-decompos}

It is well-known that for a compact $X$, the representation
$(\Pi_X, G, L^2(X))$ decomposes into a direct (infinite) sum of
irreducible representations of $G$ with finite multiplicities (see
\cite{GGPS}). We will fix one such decomposition and call it the
automorphic spectrum of $X$. We can write
\begin{equation*}\label{spec-L2X-1}
L^2(X)=(\oplus_i L_i)\oplus(\oplus_\kappa L_\kappa)\ ,
\end{equation*}
where $L_i$ are irreducible representations corresponding to Maass
forms (including the trivial representation), and $L_\kappa$ are
irreducible representations of discrete series.

For us it will be convenient to write this decomposition as the
following decomposition of the Hermitian form $P_X$ on $C^\8(X)$
\begin{equation}\label{spec-L2X}
P_X=\sum_iP_i\ +\ \sum_\kappa P_\kappa\ ,
\end{equation} where $P_i=pr_i^*(P_X)$ and
$P_\kappa=pr_\kappa^*(P_X)$.

\section{Triple products}\label{3-prod}
We introduce now our main object of study.
\subsection{Automorphic triple products}\label{3prod-def}\label{aut3prod}
Suppose we are given three automorphic representations $(\pi_j,
V_j, \nu_j)$, $j=1,2,3\ $ of $G$
\begin{eqnarray*}
\nu_j:V_j\to C^\8(X)\ \ .
\end{eqnarray*}

We define the $G$-invariant trilinear form
$l^{aut}_{V_1,V_2,V_3}:V_1\otimes V_2\otimes V_3\to\bc$
by the formula
\begin{eqnarray*}\label{aut3prod-def}
l^{aut}_{V_1,V_2,V_3}(v_1\otimes v_2\otimes v_3)= \int_{
X}\phi_{v_1}(x)\phi_{v_2}(x)\phi_{v_3}(x)d\mu_X\ ,
\end{eqnarray*}
where $\phi_{v_j}=\nu_j(v_j)\in C^\8(X)$ for any $v_j\in V_j$.

Let $(\pi, V,\nu)$ and $(\pi', V',\nu')$ be two {\it fixed} automorphic
representations of class one. For any automorphic representation
$(\pi_i,V_{\lm_i},\nu_i)$ of class one, we have the automorphic
trilinear functional
\begin{eqnarray*}
l^{aut}_{V,V',V_{\lm_i}}: V \otimes V' \otimes V_{\lm_i} \to\bc\ .
\end{eqnarray*}

In particular, the triple periods $c_i$ in (\ref{ci}) can be
expressed in terms of this form as
\begin{eqnarray}\label{c-l}
c_i= l^{aut}_{V,V',V_{\lm_i}}(e\otimes e'\otimes e_{\lm_i})\ ,
\end{eqnarray}
where $e\in V$, $e'\in V'$, $e_{\lm_i}\in V_{\lm_i}$, are standard $K$-fixed
unit vectors.

\subsection{Uniqueness of triple products} The central fact about
invariant trilinear functionals is the following uniqueness
result:

\begin{thm}{ubi} Let $(\pi_j,V_j),$ where $j=1,2,3,$
be three irreducible smooth admissible
 representations of $G$. Then
$\dim\Hom_G(V_1\otimes V_2\otimes V_3,\bc)\leq 1$.
\end{thm}

\begin{rem}{rem-uni}The uniqueness statement was proven
by A. Oksak in \cite{O} for the group $\SLC$ and the proof could
be adopted for $\PGLR$ as well (see also \cite{Mo} and \cite{Lo}).
For the $p$-adic $GL(2)$, more refined results were obtained by D.
Prasad (see \cite{Pr}). He also proved the uniqueness when at
least one representation is a discrete series representation of
$\GLR$.

There is no uniqueness of trilinear functionals for
representations of $\SLR$ (the space is two-dimensional). This is
the reason why we prefer to work with $\PGLR$.

We note, however, that the absence of uniqueness does not pose any
serious problem for the method we present. All what is really
needed for our method is the fact that the space of invariant
functionals is finite dimensional.
\end{rem}

\subsection{Model triple products}\label{mod3prod-def}
In Section \ref{modfunc}, we use an explicit model for
representations $(\pi,V)$, $(\pi',V')$ and $(\pi_{i},V_i)$ to
construct a model invariant trilinear functional. The model
functional will be given by an explicit formula. We call it the
{\it model triple product} and denote it by
$l^{mod}_{V,V',V_{\lm_i}}$, or simply $l^{mod}_{\lm_i}$, if $\pi$
and $\pi'$ are fixed.

These model functionals are defined for any three
irreducible unitary representation of principal series of $G$,
even if these are not
automorphic.

By the uniqueness principle for representations $\pi,\pi',\pi_i$,
there exists a constant $a_i=a_{V,V',V_{\lm_i}}$ such that:
\begin{eqnarray}\label{coef-a-def}
l^{aut}_{V,V',V_i}=a_i \cdot l^{mod}_{V,V',V_{\lm_i}}\ .
\end{eqnarray}

The constant $a_i$ depends on the automorphic realization of
abstract representations $\pi,\ \pi'$ and $\pi_{\lm_i}$, and on
the choice of the model functional
$l^{mod}_{\lm_i}=l^{mod}_{V,V',V_{\lm_i}}$.

From now on we will work with the coefficients
$d_i=|a_i|^2$.

\subsubsection{Exponential decay}\label{expdecay}
Relations \eqref{c-l} and \eqref{coef-a-def} give rise to a
formula for the triple product coefficients $c_i$
\begin{eqnarray*}\label{aut-mod-vect}
c_i=l^{aut}_{\lm_i}(e\otimes e'\otimes e_{\lm_i})=a_i \cdot
l^{mod}_{\lm_i}(e\otimes e'\otimes e_{\lm_i}) \ .
\end{eqnarray*}
Let us explain how one can deduce the exponential decay for the
coefficients $c_i$ using this identity.

The value of the model triple product functional
$l^{mod}_{\lm_i}(e\otimes e'\otimes e_{\lm_i})$ constructed in
Section~\ref{modfunc} is given by an explicit integral. In
\cite{BR3}, Appendix A, we evaluated this integral in terms of the
standard Euler $\G$-function by a direct computation in the model
and showed that $|l^{mod}_{\lm}(e_\tau\otimes e_{\tau'}\otimes
e_{\lm})|^2=1/\g(\lm)$, where $\g(\lm)$ is as in Section~\ref{results}.
After applying the Stirling formula to that
expression, one sees that it has an exponential decay in $|\lm|$.
Hence, in order to obtain bounds on the coefficients $c_i$, one
needs to bound coefficients $d_i=|a_i|^2$. In \cite{BR3} we showed  that the
coefficients $d_i$ are at most polynomial. This explains the
exponential decay of coefficients $c_i$. We note that the
coefficients $d_i$ encode deep arithmetic information, e.g., special
values of $L$-functions.

\section{Hermitian forms}\label{herm-forms}

\subsection{Hermitian forms and trilinear
coefficients $d_i$} \label{diag-aut} We explain now how to obtain
bounds for the coefficients $d_i$

Our method is based on the fact that these coefficients appear in
the spectral decomposition of some geometrically defined Hermitian
form on the space $E$ which is essentially the tensor product of
spaces $V$ and $V'$. This form plays a crucial role in what
follows.

More precisely, denote by $L$ and $L'$ the Hilbert completions
of spaces $V$ and $V'$, consider the unitary representation $(\Pi,
G \times G, L \otimes L')$ of the group $G \times G$ and denote by
$E$ its smooth part; so $E$ is a smooth completion of $V \otimes
V'$.

Denote by $\CH(E)$ the (real) vector space of continuous
Hermitian forms on $E$ and by $\CH^+(E)$ the cone of nonnegative
Hermitian forms.

We will describe several classes of Hermitian forms on $E$; some
of them have spectral description, others are described
geometrically.

Let $W$ be a smooth unitary admissible representation of $G$. Any
$G$-invariant functional $l: V\otimes V'\otimes W \to \bc$ defines a
$G$-intertwining morphism $T^l: V \otimes V' \to W^*$ which extends
to a $G$-morphism
\begin{eqnarray*}\label{T-W-def}
T^l: E \to \overline W \ ,
\end{eqnarray*}
where we have identified the complex conjugate space $\overline W$
with the smooth part of the space $W^*$ (see Section
\ref{complex-c}).

The standard Hermitian form (scalar product) $P = P_W$ on the space
$W$ induces the Hermitian form $\bar P$ on $\overline W$. Using
the operator $T^l$ we define the Hermitian form $H^l$ on the space
$E$ by $H^l = (T^l)^{\ast} (\bar P)$, i.e., $H^l(u)= \bar
P(T^l(u))$ for any $u\in E$.

\begin{rem}{}We note that if the representation of $G$ in the space $W$ is
irreducible and $l\not=0$, then starting with the Hermitian form
$H^l$, we can reconstruct the space $W$, the functional $l$, and
the morphism $T^l$ uniquely up to an isomorphism.
\end{rem}

\subsubsection{Forms $H_\lm$} Let us introduce a special notation for the particular case we are
interested in. For any $\lm \in i \br\cup(0,1)$, consider the
class one representation $W = V_\lm$, choose the model trilinear
functional $l^{mod}_\lm: V \otimes V' \otimes V_\lm \to \bc$
described in Section \ref{modfunc} and denote the corresponding
Hermitian form on $E$ by $H^{mod}_\lm$ or simply by $H_\lm$.
Accordingly, let $H^{aut}_i$ be the form corresponding to the
automorphic functional. We have
$H^{aut}_{i}=d_i \cdot H^{mod}_{\lm_i}$,
where $d_i=|a_i|^2=|a_{V,V',V_i}|^2$ are as in (\ref{coef-a-def}).
This is the definition of the coefficients $d_i$ we are going to
work with.

\subsection{Diagonal form $H_\Dl$} \label{diagonal}

Consider the space $C^\8(X\times X)$. The diagonal $\Dl:X \to
X\times X$ gives rise to the restriction morphism $r_\Dl :
C^\8(X\times X)\to C^\8(X)$. We define a nonnegative Hermitian
form $H_\Dl$ on $C^\8(X\times X)$ by $H_\Dl =
(r_\Dl)^{\ast}(P_X)$, i.e.,
$$H_\Dl(u)= P_X(r_\Dl(u)) = \int_{X}|r_\Dl(u)|^2d\mu_X $$
 for any $u\in C^\8(X\times X)$.

We say that $H_\Dl$ is the {\it diagonal form}.

We now consider the spectral decomposition of the Hermitian for
$H_\Dl$ (for a detailed discussion, see \cite{BR3}). Using the
spectral decomposition \eqref{spec-L2X} $P_X=\sum_iP_i\ +\
\sum_\kappa P_\kappa$ we can write $H_\Dl=\sum_i H^{aut}_i\ +\
\sum_\kappa H^{aut}_\kappa$. We have seen before that $H_i^{aut}=d_i
H_{\lm_i}$. Hence we have the following spectral identity (which is
a version of the Parseval identity)
\begin{eqnarray*}
H_\Dl=\sum_{i}d_i H_{\lm_i} + \sum_{\kappa}
H^{aut}_{\kappa} \ .
\end{eqnarray*}

Here the summation on the right is over {\it all} irreducible
unitary automorphic representations appearing in the decomposition
of $L^2(X)$ (see \eqref{spec-L2X}). The first sum is over the class
one automorphic representations (including the trivial one) and the
second sum is over the discrete series automorphic representations.

\begin{rem}{K-bessel-eq}For  most of the proof we will need
just the inequality (the Bessel inequality)
\begin{eqnarray}\label{spec-Q-Del-aut-1}
\sum_{i}d_i H_{\lm_i} \leq H_\Dl \ .
\end{eqnarray}

In order to avoid computations with discrete series, we consider
only vectors $w\in E$ which are $\Dl K$-invariant under the
natural diagonal action of $\Dl G\subset G\times G$ on $E$. For
such vectors, the inequality \eqref{spec-Q-Del-aut-1} becomes the
equality
\begin{eqnarray}\label{parseval-K} H_\Dl(w)=
\sum_{i}H^{aut}_{i}(w)=\sum_{i}d_iH_{\lm_i}(w)\ .
\end{eqnarray}
Here the
summation is over all automorphic representations of class one.

This follows from the simple fact that for a $\Dl K$-invariant
vector $w\in E$, the restriction onto the diagonal $\Dl X$ of the
automorphic realization $\nu\otimes\bar\nu( w)$ is a $K$-invariant
function on $X$, and hence orthogonal to discrete series
representations appearing in $L^2(X)$.

\end{rem}

\section{$L^4$-norm of $K$-types}

In this section we prove Theorem \ref{thm2}. We assume, for
simplicity, that the representation $V$ is a representation of the
principal series.

\subsection{Complex conjugate representation}\label{complex-c}
Our proof of Theorem \ref{thm2} is spectral, it is based on the basic
spectral identity \eqref{besseli}
 applied to the case when the representation $V'$ coincides
 with the complex
 conjugate $\bar{V}$ of the representation $V$.

We recall that for any complex vector space $V$ we can define
 the complex conjugate space ${\bar V}$. By definition, $\bar V$ is the same
 real vector space as $V$, i.e., we have a canonical bijection $c: V
 \to \bar{V}$, and the structure of the complex vector space is given by
 $\lm c(v) = c(\bar{\lm } v)$, $\lm\in\bc$. In particular, $c$ is
 an antilinear bijection.

The complex conjugate representation $(\bar{\pi}, G, \bar{V})$ naturally corresponds to any
representation $(\pi, G, V)$; unitary structure on $V$ defines a unitary structure on $\bar{V}$.

Let us note that for $\tau\in i\br$, the representation $\bar{V_\tau}$ is
 {\it canonically} isomorphic to the representation
 $V_{\bar{\tau}}$ when we consider them as spaces of functions on
$\br^2 \sm 0$ (see Section~\ref{rep-models}). The isomorphism is given
by the complex conjugation $c(v)=\bar v$.

An Hermitian
form on a space $V$ gives rise to the morphism $V\to V^+$, where
$V^+:=\overline{(V^*)}$ is the complex conjugate of the dual space.

\subsection{Complex conjugate representation in automorphic
picture}

Suppose now that we fixed an automorphic structure $\nu: V\to C^\8(X)$
 on the representation $V$. Then it
 defines the canonical automorphic structure $\bar{\nu}: \bar{V} \to
 C^\8(X)$ on the complex conjugate representation
 by the formula $\bar{\nu}(c(v)) = \overline{\nu(v)}$.

We will consider the representation $E=V\otimes\bar V$ of
 the group $G \times G$ and denote by $\nu_E = E\to C^\8(X\times
X)$ the corresponding automorphic structure on $E$ (here $\nu_E =
\nu\otimes\bar\nu$). We have the following basic claim
(compare with \ref{g-bound-proof}).

\begin{claim}{aut-positive}For any vector $v
\in V$ consider the vector $w = v \otimes \bar{v}=v\otimes c(v) \in E$ and the
corresponding function $\nu_E(w)$ on $X \times X$. Then the
restriction $u = r_\Dl(\nu_E(w))$ of this function to the diagonal
$\Dl X$ is a non-negative function on $X$, and $H_\Dl(w) =
||u||^2_{L^2(X)} = ||\nu(v)||^4_{L^4(X)}$.
\end{claim}

This follows from the observation that $u(x) = \nu(v)(x) \cdot \overline
{\nu(v)(x)} = |\nu(v)(x)|^2$.
\qed

\subsection{$K$-types} We assume that $V = V_\tau$ is a representation of
the principal series. All the necessary computations will be done in
the circle model $V_\tau\simeq C^\8_{even}(S^1)$ (i.e., we realize a
vector in $V$ as a smooth function $f$ of the angular parameter
$t\in\br$ such that $f(t+\pi)=f(t)$). The invariant unitary
Hermitian form on $V$ is given by
$||f||^2=\frac{1}{\pi}\int_0^\pi|f(t)|^2dt$.

Let $e_n=\exp(int)$, where $n\in 2\bz$, be an orthonormal basis of
$K$-types in the space $V_\tau$ (all weights are even since we work
with the group $G=\PGLR$).

Consider the space $\bar V_\tau$. We have a natural identification
$\bar V_\tau\simeq V_{-\tau}$ induced by the realization of these
spaces as spaces of functions on $\br^2\setminus 0$.

We denote by $\{e'_n=\overline{ e_{-n}}\}_{n\in 2\bz}$ the
corresponding complex conjugate basis for $\bar V_\tau\simeq
V_{-\tau}$. Under the natural identification $V_{-\tau}\simeq
C^\8_{even}(S^1) $, we have $e'_n=\exp(int)$ as before.

\subsection{Test vectors}\label{sect-test} In the Introduction
(see formula \eqref{test-vectors}) we defined two families of test
vectors central for our proof of the subconvexity. We repeat this
construction.

For any $n\in 2\bz$, $n\geq 0$, we consider two vectors in
$E=E_\tau=V_\tau\otimes V_{-\tau}$ given by
\begin{eqnarray*}\label{w1}
w_n=e_n\otimes e'_{-n}\ ,\  {\rm and}\ \  \tilde{w}_n=w_n+w_{n+2}\ .
\end{eqnarray*}

We note that in the model $V_\tau\otimes
V_{-\tau}\simeq C^\8_{even,even}(S^1\times S^1)$ these vectors are
represented by the functions $w_n(x,y)=e^{in(x-y)}$ and $\tilde
w_n(x,y)=(1+e^{i2(x-y)})e^{in(x-y)}$.

In Section \ref{g-bound-proof} we have proven the basic geometric
bound \eqref{(star)} for these vectors
 $$H_\Dl(w_n)\leq H_\Dl(\tilde
w_n)\ .\eqno{(\star)}$$

\subsection{Main Proposition} Our main claim is the following

\begin{prop}{main-prop}There exists a positive constant $D$ such that
$$H_\Dl(\tilde w_n)\leq D\ , \eqno{(\natural)}$$
for all $n$.
\end{prop}\\
We prove this proposition in Section \ref{proof-K-types-4norm}.

\begin{rem}{}The bound ($\star$) is of
a geometric nature as it concerns the form $H_\Dl$ defined on the
automorphic space $X$ and appeals to the automorphic realization
of $V$ in $C^\8(X)$. On the other hand, our proof  of the bound
$(\natural)$ is purely {\it spectral}, despite its geometric
appearance.
\end{rem}

\subsection{Proof of Theorem \ref{thm2}} Proposition \ref{main-prop} and the geometric bound $
H_\Dl(w_n)\leq H_\Dl(\tilde w_n)$ (see \eqref{geom-bound}) imply the
bound in Theorem \ref{thm2} for $L^4$-norm of $K$-types. Namely,
from Claim \ref{aut-positive} we see that
\begin{eqnarray}\label{L-4-app}||\nu(e_n)||_{L^4(X)}^4=H_\Dl(w_n)
\leq H_\Dl(\tilde w_n)\leq D\ ,\end{eqnarray} for some $D$
independent of $n$. \qed

\begin{rem}{D-depend}The method presented in this paper allows one to  give an effective estimate
for the constant $D$ in Theorem \ref{thm2} (and in Proposition
\ref{main-prop}). It depends  on geometry of the Riemann surface $Y=X/K$  and on the parameter
$\tau$ of the representation $V$. Namely,  the following
bound $$D\leq C\cdot\frac{{\rm vol}(Y)}{{\rm vol}(B_{i(Y)})}\cdot (1+|\tau|)^2\ ,$$ should hold for some absolute
constant $C>0$. Here $B_{i(Y)}$ is a hyperbolic ball of the radius equal to the injectivity radius $i(Y)$ of $Y$.

Careful execution of ideas presented in \cite{Re2} should give a better bound $||\nu(e_n)||_{L^4}^4\leq A_Y\cdot\left(1+\frac{1+|\tau|^{3/2}}{1+|\tau|+|n|}\right)
$, with a constant $A_Y$ depending on $Y$ only.  In particular, from this would follow that for $|n|\gg |\tau|^{3/2}$, $L^4$-norm of a $K$-type $e_n\in V_\tau$ is uniformly bounded {\it independently} of $\tau$. For $n=0$ (i.e., for the Maass form $\phi_\tau$ on $Y$),  such a bound is consistent with the general PDE bound of C. Sogge \cite{So} (i.e., $||\phi_\tau||^4_{L^4}\leq C''|\tau|^\haf$). One expects that the correct bound is $||\nu(e_n)||_{L^4}^4\ll (1+|n|+|\tau|)^\eps$ for any $\eps>0$. For  a congruence subgroup $\G$ and Hecke-Maass forms (i.e., $n=0$), this is the result of P. Sarnak and T. Watson (unpublished). We plan to discuss these issues elsewhere.

\end{rem}

\section{Proof of Proposition \ref{main-prop}}

\subsection{Spectral Lemma}\label{proof-geom-bd} Our proof is based on
 the following spectral
bounds (these are bounds
 \eqref{(ii)} from the Introduction).

Recall that  we set $\mathcal{S}=2(|\tau|+|\tau'|)+1$ (in fact in this section we can assume that $\tau'=-\tau$).

\begin{lem}{lem-II}There exists a constant
$C$ such that for any $n\in 2\bz$, the following spectral
bounds hold
\begin{description}
 \item[$(II_1)$] $H_\lm(\tilde w_n)\leq
C\cdot(1+ |n|)\inv|\lm|^{-1}+C|\lm|^{-3}$ for all $\lm$ satisfying
$\mathcal{S}\leq|\lm|\leq 4|n|$,
 \item[$(II_2)$] $H_\lm(\tilde w_n)\leq C|\lm|^{-3}$ for all $\lm$
 satisfying $|\lm|\geq
 4|n|$.
\end{description}
\end{lem}

The model Hermitian forms $H_\lm$ on $E$ are defined explicitly for
every $\lm\in i\br$ as in Section \ref{modfunc}. The proof of the
lemma amounts to a routine application of the stationary phase
method and the van der Corput lemma (see Section \ref{proof-bounds-II-1-2} ).
In fact, the restriction $|\lm|\geq \mathcal{S}$ is purely technical. One can obtain good bounds for the value of $H_\lm(\tilde w_n)$ for all $\lm$. We will not need this in what follows. The constant $C$ in the lemma above satisfies a bound $C\leq C'\cdot \mathcal{S}$ for some absolute constant $C'$.

\subsection{Proof of Proposition \ref{main-prop}}\label{proof-K-types-4norm}
For any given $n$, the function $\nu_E (\tilde w_n)$ is a bounded
smooth function on $X\times X$ and hence $H_\Dl(\tilde w_n)$ is
well-defined. We have to show  that it is bounded by some constant
$D$ independent of $n$.

As could be seen from our construction in Section \ref{sect-test},
vectors $\tilde w_n$ are $\Dl K$-invariant. It follows from the
discussion in Remark \ref{K-bessel-eq} that for such vectors, we
have the following Parseval identity \eqref{parseval-K}
$$H_\Dl(\tilde w_n)=\sum_{i}H^{aut}_{i}(\tilde
w_n)=\sum_{i}d_iH_{\lm_i}(\tilde w_n)\ .$$

Here the sum is over the spherical spectrum $I =\{\lm_0,\lm_1,...\}$.
Let $k_0\in\mathbb{N}$ be such that $2^{k_0}\leq \mathcal{S}\leq 2^{k_0+1}$. We decompose the spherical spectrum $I$ as a union of subsets $I_{k_0},\
I_{k_0+1},\ ...$ (dyadic intervals) according to the absolute value of
$|\lm|$, and estimate the contribution of each of these subsets.

Namely, we consider dyadic subsets $I_k$ of the spectrum $I$ defined by
$I_{k_0} = \{ \lm \in I\  | \ |\lm| < 2^{k_0+1} \}$ and $I_k = \{ \lm \in I\  | \
|\lm| \in [2^k, 2^{k+1}) \}$ for $k > k_0$.

Notice that all exceptional spectra that correspond to representations
of the complementary series and to the trivial representation is contained
in the interval $I_{k_0}$ (we call it the low spectrum).
All the other intervals contain only imaginary values of $\lm$ which
correspond to representations of the principal series.

We have $H_\Dl(\tilde w_n) = \sum_{k\geq k_0} H_k$,
where $H_k = \sum_{\lm_i \in I_k} d_iH_{\lm_i}(\tilde w_n)\ .$

\subsubsection {Estimate of $H_k$ for $k >k_0$}
The idea of the proof is that on the interval $I_k$ the function
$H_\lm(\tilde w_n)$ is more or less constant, so we will not lose
much when we replace it by its maximal value.

According to the bound $(II_{2})$, Lemma \ref{lem-II}, we see that
for $\lm \in I_k$ we have a bound $H_\lm(\tilde w_n) \leq M_k$ where
$M_k = C (n^{-1} 2^{-k} + 2^{-3k})$ for $k$ satisfying $2^k < 4n$,
and $M_k = C 2^{-3k}$ for $k$ satisfying $2^k \geq 4n$. Here $C$ is
a universal constant that depends only on $\tau$.

According to the mean-value bound \eqref{mean-value-convex} we have
$\sum_{\lm_i \in I_k} d_i \leq A \ 2^{2k}$. Hence we arrive at the
bound $H_k \leq 2^k A M_k$. This implies that
$$\sum_{k>k_0 } H_k \leq AC\left( \sum_{k>0} 2^{-k} + \sum_{2^k < 4n} 2^k n^{-1}\right)
\leq AC (1 + 8)\ .$$

\subsubsection{ Estimate of the low spectrum contribution $H_{k_0}$}\label{low-spec}
 We claim that the sum $H_{k_0}$ is bounded by some constant $D'$ which depends
 only on the geometry of the space $Y$. In principle we could apply
 to this case the spectral argument similar to the one described above.
 However this would lead to some unpleasant computations with the exceptional
 spectrum.
For that reason we prefer to give  the following more geometric
argument.

The vector $\tilde w_n\in E$ is a $\Dl K$-invariant vector. Hence the
corresponding function $b=\nu_E(\tilde w_n)|_{\Dl X}=
|\phi_n|^2+|\phi_{n+2}|^2$ is a $K$-invariant function and we can
view it as a function on $Y$. Moreover, we can compute its $L^1$-norm
on $Y$
\begin{eqnarray*}
||b||_{L^1(Y)}=\int_Y\left(|\phi_n|^2+|\phi_{n+2}|^2\right)dv= 2\ .
\end{eqnarray*}

Consider the subspace $R={\rm span}\{\phi_{\lm_i}\ |\ |\lm_i|<2^{k_0+1}\}
\subset C^\8(Y)$. This is a finite-dimensional vector space
consisting of smooth functions. Since the space $R$ is
finite-dimensional we can bound the supremum norm on this space
$||\cdot||_\8$ by $L^2$-norm, i.e., there exists a constant $C_R$
such that $||f||_\8 \leq C_R ||f||_{L^2(Y)}$ for all functions $f
\in R$.

\begin{claim}{} $H_{k_0} \leq 4 C_R^2$. \end{claim}

Indeed, by  definition $H_{k_0} = ||a||^2_{L^2(Y)}$, where the vector $a
\in R$ is the orthogonal projection $a = pr_R(b)$ of the vector $b$
onto the subspace $R \subset L^2(Y)$.

Thus we have \begin{eqnarray*}H_{k_0}^2 = |\langle a,a\rangle|^2 = |\langle b,a\rangle|^2 \leq  ||b||_{L^1(Y)}^2 \cdot
 ||a||_\8^2 \leq 4C_R^2 \cdot  ||a||_{L^2(Y)}^2 = 4\ C_R^2 \cdot H_{k_0}\ .\end{eqnarray*}
This implies the claim and finishes the proof of the proposition. \qed

\begin{rems}{proof-remks}1.
It is not difficult to bound the constant $C_R$ in the proof above
   in terms of the geometry of the Riemann surface $Y$ and the parameter
   $\mathcal{S}$. For example, suppose we found a number $r < 1$ that is smaller
   than the injectivity radius of $Y$. Then one can show that
$C_R^2  \leq 100 (\mathcal{S}^2 + (1/r)^2) {\rm vol}(Y)$,
where ${\rm vol}(Y)$ is volume computed with respect to the standard hyperbolic metric.

2.  The proof of  Proposition \ref{main-prop} given above could be easily extended  to the
case of a general finite co-volume lattice $\G \subset G$. In fact, the only place where we implicitly used compactness of $X$ is in the proof of the mean-value bound \eqref{mean-value-convex} which we quoted from \cite{BR3}. However, in \cite{BR1} we proved similar bound for a general finite co-volume lattice and cuspidal functions $\phi$ and $\phi'$.

For a general finite co-volume lattice, the spectral decomposition of the Laplace-Beltrami operator
on $Y = \G \backslash \uH$ is given by a collection of
eigenfunctions $\phi_z$, where the parameter $z$ runs through some
set $Z$ with the Plancherel measure $d\mu$. The spectral set $Z$
has discrete points which correspond to eigenfunctions (Maass
forms) $\phi_z\in L^2(Y)$, and the continuous part which corresponds to
eigenfunctions coming from the unitary Eisenstein series. The
collection $\{\phi_z\}_{z\in Z}$ defines a transform $\hat
u(z)=<u,\phi_z>$ for every $u\in C^\8_c(Y)$. The main property of this
transform is the Plancherel formula $||u||^2_{L^2(Y)}=\int_Z |\hat
u(z)|^2 d\mu .$

Let us fix two Maass {\it cusp forms} $\phi$ and $\phi'$ on $Y$.
For every $z \in Z$, we define the parameter $\lm_z\in\bc$ and the
coefficient $d_z$ in the same way as before. In \cite{KS} the following mean-value bound was obtained (improving on our result in \cite{BR1})
\begin{eqnarray*}
\int_{T\leq |\lm_z|\leq 2T} d_z \ d\mu \leq A\left(\ln(T)\right)^\frac{3}{2}\cdot T^{2 }
\ .
\end{eqnarray*}

The proof given in present paper, together with  the above mean-value bound,
gives the following bound for $L^4$-norm of $K$-types in a {\it class one fixed cuspidal} representation  $\nu:V\to C^\8(X)$
 \begin{eqnarray*}
 ||\nu(e_n)||_{L^4(X)}\leq D\left(\ln(2+|n|)\right)^\frac{3}{2}\  \rm{for\ all}\ n.
 \end{eqnarray*} This is our analog of Theorem \ref{thm2}  for non-uniform lattices. In particular we do not know wether  $L^4$-norm of $K$-types are uniformly bounded for a non-uniform lattice.

The bound on  $L^4$-norm of $K$-types  implies  as before that the following subconvexity bound holds for a general finite co-volume lattice
\begin{eqnarray*}
\int_{Z_T} d_z \ d\mu \leq B\left(\ln(T)\right)^\frac{3}{2}\cdot T^{5/3}\ , \ {\textrm{where}} \
 Z_T = \{z \in Z \ |\ |\lm_z| \in I_T \}
\ ,\end{eqnarray*} for some constant $B>0$.
\end{rems}{}
\\

The rest of the paper is devoted to the proof of spectral bounds
 $(II_{1,2})$ from Lemma \ref{proof-geom-bd} and the lower bound \eqref{lower-bd}. This will be done using
computations in the explicit model of irreducible representations.
As a preparation we start with an explicit construction of model
Hermitian forms $H_{\lm}$.


\section{Model trilinear functionals}
\label{sect2}

\subsection{Model trilinear functionals} \label{modfunc}

In this section we briefly recall our construction from \cite{BR3}
of model trilinear invariant functionals.

For every $\lm \in \bc$, we denote by $(\pi_\lm,V_\lm)$ the smooth
class one representation of the generalized principle series of the
group $G=\PGLR$ described in Section \ref{rep-models}. As a vector space $V_\lm$
is isomorphic to the space of smooth even functions $C^\8_{even}(S^1)$ on $S^1$.

We describe the {\it model} invariant trilinear functional using
this geometric model. Namely, for three given complex numbers
$\tau,\ \tau',\ \lm$, we explicitly construct a nontrivial trilinear
functional $\ l^{mod}: V_{\tau} \otimes V_{\tau'} \otimes V_{\lm}
\to \bc$ by means of its kernel. In the circle model, the
trilinear functional on the triple $V_\tau,\ V_{\tau'},\ V_\lm$ is
given by the following integral:
\begin{equation*}\label{circleintegral}
l^{mod}_{\pi,\pi',\pi_\lm}(f_1\otimes f_2\otimes f_3)= (2
\pi)^{-3} \int\limits_{(S^1)^3}
f_1(x)f_2(y)f_3(z)K_{\tau,\tau',\lm} (x,y,z) dx dy dz\ ,
\end{equation*}
with the kernel
\begin{multline}\label{kernel-S}
K_{{\tau,\tau',\lm}}(x,y,z)=
|\sin(x-y)|^{\frac{-\tau-\tau'+\lm-1}{2}}|\sin(x-z)|^{\frac{-\tau+\tau'-\lm-1}{2}}
|\sin(y-z)|^{\frac{\tau-\tau'-\lm-1}{2}}.
\end{multline}

Here $x, y , z$ are the standard angular parameters on the circle $S^1$.
As we verified in \cite{BR3} this defines a non-zero $G$-invariant
functional.

\begin{rem}{rem-trip}1. In general the integral
defining the trilinear functional is often divergent and the
functional should be defined using regularization of this integral.
There are standard procedures how to make such a regularization (see
\cite{G1}). Fortunately, in the case of class one unitary
representations, all integrals converge absolutely, so we will not
discuss the regularization procedure.

2. We do not have a similar simple formula for the trilinear
invariant functional when at least one representation is a
representation of discrete series. This is because we do not know a
simple ``geometric" model for representations of discrete series. As
a result it is more cumbersome to carry out explicit computations in
that case. Another problem we have to face is that the results of
\cite{BR3} have not been extended yet to cover the discrete series.

Nevertheless, we expect our methods to carry out for discrete
series as well and to produce corresponding subconvexity bounds, and bound for $L^4$-norms of $K$-types.
\end{rem}

\subsection{Reduction for $\Dl K$-invariant vectors}\label{K-reduction}
In what follows, we only need to deal with $\Dl K$-invariant vectors
in $E\simeq V_\tau\otimes V_{\tau'}$.  For such vectors, we can reduce
the integral \eqref{kernel-S} representing the model
invariant functional, and hence the Hermitian form $H_\lm$ to the integral in one variable.

Namely,
let $l_\lm^{mod}:E\otimes V_\lm\to\bc$ be the model trilinear
functional introduced in Section~\ref{modfunc},
$T_\lm=T_\lm^{mod}:E\to V_{-\lm}$ be the corresponding map, and
$H_\lm$ the model Hermitian form on $E$ obtained from the
composition of $T_\lm$ with the invariant unitary form on
$V_{-\lm}$. We assume that $V_\lm$ is a representation of the
principal series since we are only interested in the case when
$|\lm| \geq \mathcal{S}$. In this case, the unitary form on $V_{\lm}\simeq
C^\8_{even}(S^1)$ is the standard normalized unitary form on
$L^2(S^1)$.

Let $w\in E\simeq C^\8_{even,even}(S^1\times S^1)$ be a $\Dl
K$-invariant vector. Since it is $\Dl K$-invariant it can be
represented by a function of one variable
 $c=x-y$: $w(x,y)=u(c)$, where $ u\in C^\8_{even}(S^1)$.
We claim that the estimate of $H_\lm(w)$ could be reduced to
an estimate of an integral in one variable. Namely, on the space of $\Delta K$-invariant vectors in
$E$ the form $H_\lm$ has rank $1$, i.e., it is equal to the absolute value squared of
some functional $b_\lm$ on $C^\8(S^1)$. More precisely, we have the following

\begin{lem}{H-to-b}Fix $\tau,\ \tau'\in i\br$ as before and assume that
$\lm\in i\br$.
There exists an $L^1$ function $l_\lm$ on $S^1$ such
  that for any function $u\in C^\8_{even}(S^1)$ and for the corresponding
  vector $w(x,y)=u(x-y) \in E$, we have $H_\lm(w)= |b_\lm(u)|^2$ where
   $b_\lm (u) =
\int  l_\lm (c) u(c) dc$.
\end{lem}

\Proof Since the vector $w$ is $\Delta K$-invariant its image
$T_{\lm}(w)\in V_\lm$ is proportional to the standard unit
$K$-invariant vector $e_\lm$. The proportionality coefficient
$b_\lm(u)$ equals
\begin{eqnarray*}T_\lm(w)(0) = (1/ 2\pi)^2 \int K_{\tau,\tau',\lm}(x,y,0) w(x,y) dx dy =
 1/2\pi  \int_{S^1}  l_\lm(c) u(c) dc\ ,
 \end{eqnarray*}
 where
 \begin{eqnarray}\label{def-l-lm} l_\lm(c) = \frac{1}{2\pi} \int_{S^1} K_\lm(y + c, y, 0) dy \end{eqnarray}
and  $K_\lm(x,y,z)$ is the kernel of the model trilinear
functional defined in \eqref{kernel-S}.

 Thus we see that
 $H_\lm(w) = ||T_\lm(w)||^2 = |b_\lm(u)|^2$.\Endproof
\begin{rem}{}Uniqueness of trilinear functionals implies that $b_{-\lm}=a(\lm)\cdot b_\lm$ for some scalar $a(\lm)\in\bc^\times$. It is also easy to see that $|a(\lm)|=1$.
\end{rem}


\section{Proof of spectral bounds}

\subsubsection{A convention} In what follows we will study asymptotic
behavior for various oscillating integrals. We will consider
expansions consisting of a main term and a remainder. We will bound
corresponding remainders in terms of $C^N$-norms.

We will use the following notations.  We consider functionals on $C^\8(\br)$ of the form $I_\Lambda (\phi)=\int _\br k_\Lambda( x)\phi(x)dx$ where $\phi\in C^\8(\br)$ (usually with compact support). Here $k_\Lambda( x)\in L^1(\br)$ is a kernel function depending on a set of parameters $\Lambda\in\br^n$.
We consider approximations of such functionals of the form
$I_\Lambda(  \phi)=I^0_\Lambda(\phi)+RI_\Lambda( \phi)$ where we call $I^0_\Lambda( \phi)$ the main term and
$RI_\Lambda( \phi)$ the remainder. Usually, the main term will be  given by the stationary phase method (i.e., it will be given by a functional which is a weighted sum of $\dl$-functions at points corresponding to critical points of the phase of $k_\Lambda$). We will consider bounds for $RI_\Lambda(\phi)$ in terms for $C^N$-norms of function $\phi$.
For $\phi\in C^\8(a,b)$ and an integer $N \geq 0$, we will  denote by $||\phi||_{C^N}$
the $C^N$-norm of $\phi$ defined by $||\phi||_{C^N}=\sup\limits_{0\leq
m\leq N,\ x\in (a,b)}|\phi^{(m)}(x)|$.

\subsection  {Estimate of the functional $b_\lm$.}
In Section \ref{K-reduction} we have reduced estimates of the form $H_\lm$ to the estimates of
the functional $b_\lm$. We will be interested in the case when the function $u$ from Lemma \ref{H-to-b} has a form
$u(c) = \phi(c) e^{inc}$, where $\phi$ is  a fixed smooth function
and $n \in 2 \bz$ is a parameter. We can consider the expression
$b_\lm(u)$ as a functional $F_{\lm,n}$  on the space $C^{\8}(S^1)$
which depends on two parameters $\lm$ and $n$. This functional is given by
\begin{eqnarray}\label{F_lm-n} F_{\lm,n}(\phi) : = \int_{S^1} l_\lm (c) e^{inc} \phi(c) dc\ .\end{eqnarray}

The main technical difficulty in evaluating this functional is that we have to give estimates for the values of this functional that are uniform in two parameters $\lm$ and $n$.

Recall that we set $\mathcal{S}=2(|\tau|+|\tau'|)+1$ and assume that $|\lm|\geq \mathcal{S}$.  Using the symmetry of functional $F_{\lm,n}$, we will show  that it is enough to consider the case  when $n \in 2 \bz_+$ and $\lm = it,\ t \geq \mathcal{S}$. It turns out that under these conditions the functional $F_{\lm,n}$ is almost proportional to an elementary functional $\phi \mapsto \phi(c_0)$  where $c_0 = \pi /2$.

We have the following

\begin{prop}{H-to-Airy} Consider the functional $F_{\lm,n}$ when $n \in 2 \bz_+$ and $\lm = it,\ t \geq \mathcal{S}$. We  have the following estimates of the values of this functional in terms of $C^N$-norms on $C^\8(S^1)$. There exists $C>0$ such that
\begin{enumerate}
  \item  If $t \geq 4n$ then
        $|F_{\lm,n}(\phi)| \leq C ||\phi||_{C^3}\cdot t ^{-\frac{3}{2}}$.
  \item If $t < 4n$ we have an approximation $F_{\lm,n}(\phi)=F^0_{\lm,n}(\phi)
       + RF_{\lm,n}(\phi)\ ,$ where the main term is given by
$F^0_{\lm,n}(\phi)=A(\lm,n) \phi(c_0)$, and the error term satisfies the bound
\begin{eqnarray*}|RF_{\lm,n}(\phi)| \leq C ||\phi||_{C^2}\cdot t^{-\haf}(1+ n)^{-\haf} + C ||\phi||_{C^3}\cdot t ^{-\frac{3}{2}}\ .
\end{eqnarray*} \end{enumerate}
The coefficient $A(\lm,n)$ is given by  $A(it,n) = t^{-\frac{5}{6}}\mathbb{A}(t^{-\frac{1}{3}}(2n -t))$, where  $\mathbb A$ is  the classical Airy function (see \cite{Mag}, \cite[Section 7.6]{He}).
\end{prop}

We will prove this  proposition in Section \ref{proof-prop-H-to-Airy} by carefully estimating the oscillating integral defining the functional $F_{\lm,n}(\phi)$. For the constant $C$ above we can obtain a bound of the form $C\leq C' \mathcal{S}$ for some absolute constant $C'$.

\subsection{Proof of the spectral bound \eqref{lower-bd}}\label{proof-i}
We repeat the construction of the test vector $w_T$ in \eqref{test-vectors1}.
We assume that $V = V_{\tau}$, $V' =
V_{\tau'}$ for some $ \tau,\ \tau' \in i \mathbb{R}$.
We choose an orthonormal basis $\{e_n\}_{n\in 2\bz}$ in $V$
consisting of $K$-types and similarly
an orthonormal basis $\{e'_n\}$ in $V'$.

For a given $T\geq \mathcal{S}$, we choose even $n\geq 0$ such that $|T-2n|\leq 10$,
and set $ w_T:=e_n\otimes e'_{-n}$.

Using the reduction from Section \ref{K-reduction}, we see that the vector $w = w_T$ corresponds to a function $u(c) = e^{inc}$.
Hence we have $H_\lm(w) = |F_{\lm,n}(\phi)|^2$, where $\phi \equiv 1$.

 From (2) in Proposition \ref{H-to-Airy} we see that
$F_{\lm,n}(\phi) = A(\lm,n) \phi(c_0) + RF_{\lm,n}(\phi)$. In this case we
have $|RF_{\lm,n}(\phi)| \leq C (1+|n|)^{-1}$, $\phi(c_0) = 1$. The Airy function $\mathbb{A}$ is a smooth {\it non-vanishing} at $0$ function (\cite{Mag}, \cite[Section 7.6]{He}). Hence there are constants $b,\ c>0$ such that $|\mathbb{A}(x)|\geq c$ for all $|x|\leq b$. This implies that   $|A(it,n )| \geq ct^{-5/6}$ for $|2n-t|\leq b t^{-\frac{1}{3}}$. Hence in the approximation of $F_{\lm,n}(\phi)$ stated in  Proposition \ref{H-to-Airy}   (2),   the main term $A(it,n )\phi(c_0)$ dominates the reminder $RF_{\lm,n}(\phi)$. The lower bound \eqref{lower-bd} follows.\qed

\subsection{Proof of  Lemma \ref{lem-II}, ($II_{1,2}$)}\label{proof-bounds-II-1-2}
We assume that $V'\simeq \overline V$, i.e., $\tau=-\tau'$. Let $n\in 2\bz$ and $\lm\in i\br$, $|\lm|\geq \mathcal{S}$, and $\tilde w=\tilde w_n$ as in Section \ref{sect-test}.
As in Section \ref{proof-i},  we have $H_\lm(\tilde w) = |F_{\lm,n}(\tilde \phi)|^2$, where
$\tilde \phi(c) = 1 + e^{2ic}$. This time we are looking for a uniform in $n$ upper bound valid for
{\it all} $|\lm|\geq \mathcal{S}$.

We need to bound the integral $F_{\lm,n}(\tilde \phi)$. From the form of integral \eqref{F_lm-n} it follows
that it is enough to consider the case  $n\geq 0$ and $\rm{Im}(\lm)\geq 0$. Indeed, using the change of variables $c\mapsto-c$ in integral \eqref{F_lm-n}, we can assume that $n\geq 0$.
Considering the complex conjugate to $l_\lm$, we can assume that $\rm{Im}(\lm)\geq 0$.

Hence we can apply Proposition \ref{H-to-Airy}. We have  $\tilde \phi(c_0) = 0$,  and  hence $F_{\lm,n}(\tilde \phi) =
RF_{\lm,n}(\tilde \phi)$. Thus estimates  in Lemma \ref{lem-II} ($II_{1,2}$), directly follow from the
Proposition \ref{H-to-Airy}. \qed


\section{Proof of Proposition \ref{H-to-Airy}}\label{proof-prop-H-to-Airy}

\subsection{Proof of Proposition \ref{H-to-Airy}}\label{K-reduction-sect}
We consider  the oscillating integral
   $F_{\lm,n}(\phi)  = \int l_\lm (c) e^{inc} \phi(c) dc\ $.
   One of the difficulties in evaluating this functional is
   that its kernel function $l_\lm$ is not an elementary function.

   However, since the function $l_\lm$  itself is defined by an oscillating
   integral,  we can approximate it by an elementary function $k_\lm$
   which  is the sum of main term contributions from critical points
   of this oscillating integral.

\subsubsection{Approximation of the kernel $l_\lm$} We have the following\\
\begin{lem}{claim-H-lm}Fix $\tau,\ \tau'\in i\br$ and $\mathcal{S}$ as before and assume
that $\lm\in i\br$, $|\lm|\geq \mathcal{S}$. There exists a constant $C>0$ depending on $\tau$ and $\tau'$, such that we have the following approximation
 \begin{eqnarray} \label{H-lm-asym}
l_\lm(c)  =a_\lm\cdot |\lm|^{-\haf} k_{\lm}(c) +r_{\lm}(c)\ ,
\end{eqnarray} where $a_\lm=e^{i\frac{\pi}{4}}2^{1+\hlm}$ and the kernel
$k_{\lm}(c)$ is given by an explicit formula
$k_{\lm}(c)=A(c)m_{\lm}(c)$ with
 \begin{eqnarray} \label{H-lm-asym2}\hskip 1cm
A(c)=
|\sin(c)|^{{\frac{-\tau-\tau'-1}{2}}}\quad
,\quad m_{\lm}(c)=|\sin(c/2)|^{-\hlm}|\cos(c/2)|^\hlm\ ,
 \end{eqnarray}
and the error term $r_{\lm}(u)$ satisfies the bound
\begin{eqnarray}
|r_{\lm}(c)|\leq C|\lm|^{-\frac{3}{2}}|\sin(c)|^{-\haf}|\ln(|\sin(c/2)\cos(c/2)|)|\ .
\end{eqnarray}
\end{lem}

We will prove this lemma in Section \ref{l-k-aprox-lem-proof}. For the constant $C$ above we can obtain a bound of the form $C\leq C' \mathcal{S}$ for some absolute constant $C'$.

Using this approximation we can approximate the functional
$F_{\lm,n}$ by a simpler functional defined for $n\in 2\bz$ and $\phi\in C^\8_{even}(S^1)$, by \begin{eqnarray}\label{G-def} G_{\lm,n}(\phi): =
  \int_{S^1} k_\lm(c) e^{inc} \phi(c) dc=2\int_0^\pi k_\lm(c) e^{inc} \phi(c) dc\ .
  \end{eqnarray}

The lemma above implies \\
\begin{cor}{}There exists a constant $C'=C'(\tau,\tau')>0$ such that
\begin{eqnarray}\label{F-G}|F_{\lm,n}(\phi)- a_\lm|\lm|^{-\haf}\cdot G_{\lm,n}(\phi)|\leq C'||\phi||_{L^\8(S^1)}\cdot |\lm|^{-3/2}\ ,
\end{eqnarray} for all $|\lm|\geq \mathcal{S}$.
\end{cor}

Hence  Proposition \ref{H-to-Airy} follows from an appropriate estimate for the functional $G_{\lm,n}(\phi)$.

\subsubsection{Estimate for  $G_{\lm,n}(\phi)$} We have the following estimate for the functional  $G_{\lm,n}$ defined in \eqref{G-def}.

\begin{prop}{G-to-Airy} Consider the functional $G_{\lm,n}$ when $n \in 2 \bz_+$ and $\lm = it,\ t \geq \mathcal{S}$. There exists a constant $C>0$ depending on $\tau$ and $\tau'$, such that we have the following estimates
\begin{enumerate}
  \item  If $t \geq 4n$ then
        $|G_{\lm,n}(\phi)| \leq C||\phi||_{C^3}\cdot t ^{-3}$,
  \item If $t < 4n$ then we have an approximation $G_{\lm,n}(\phi)=G^0_{\lm,n}(\phi)
      + RG_{\lm,n}(\phi)\ ,$ where the main term is given by $G^0_{\lm,n}(\phi)=A(\lm,n) \phi(c_0)$,
and the error term $RG_{\lm,n}(\phi)$ satisfies the bound
\begin{eqnarray*}|RG_{\lm,n}(\phi)| \leq C ||\phi||_{C^2}\cdot (1+n)^{-1/2} + C ||\phi||_{C^3}\cdot t ^{-\frac{3}{2}}\ . \end{eqnarray*}\end{enumerate} The coefficient $A(\lm,n)$ is given by
 $A(it,n) = t^{-\frac{1}{3}}\mathbb{A}(t^{-\frac{1}{3}}(2n -t))$.
\end{prop}

This proposition and  bound \eqref{F-G} imply Proposition \ref{H-to-Airy}. This finishes the proof of Proposition \ref{H-to-Airy}. \qed

\subsection{Proof of Lemma \ref{claim-H-lm}}\label{l-k-aprox-lem-proof}\label{kernel-aprox-app}\label{e-kernel}
We prove  the claims in the lemma  by  essentially straightforward application of the stationary phase
method in the form explained in Appendix~\ref{red-hyp-int-sect}. In order to estimate the error of this
approximation we use the standard integration by parts argument.

 To compute the approximation $k_\lm$ of $l_\lm$, we consider for fixed
$\tau,\ \tau'\in i\br$ and for $|\lm|\geq \mathcal{S}$, $\lm\in i\br$, the integral \eqref{def-l-lm}:
\begin{eqnarray*}\label{int-red1-B}
l_\lm(c)=(2\pi)^{-\haf}\int\limits_{S^1}K_{{\tau,\tau',\lm}}(y+c,y,0)dy&&=\\
(2\pi)^{-\haf}|\sin(c)|^{\frac{-\tau-\tau'+\lm-1}{2}}\int\limits_{S^1}&&\!
|\sin(y+c)|^{\frac{-\tau+\tau'-\lm-1}{2}}
 |\sin(y)|^{\frac{\tau-\tau'-\lm-1}{2}}dy\\
&=&|\sin(c)|^{{\frac{-\tau-\tau'+\lm-1}{2}}}l'_\lm(c) \ ,
\end{eqnarray*} where the kernel $K_{{\tau,\tau',\lm}}$ is as in
\eqref{kernel-S},
and we denote by $l'_\lm$ (suppressing the dependence on $\tau,\ \tau'$) the function
\begin{eqnarray}\label{int-red2-B}
\hskip 1.5cm
l'_\lm(c) =(2/\pi)^{\haf}\int\limits_{t\in\br/\pi\bz}
|\sin(t+c/2)|^{\frac{-\tau+\tau'-\lm-1}{2}}
|\sin(t-c/2)|^{\frac{\tau-\tau'-\lm-1}{2}}dt\ .
\end{eqnarray}
To find the asymptotic of the integrals of the type of $l'_\lm(c)$ is a
problem in classical analysis. We view the integral \eqref{int-red2-B} as a
one-dimensional integral
(in $t$) with parameters $\lm$ and $c$. We treat such integrals in Appendix
\ref{red-hyp-int-sect} where we show that the main term (i.e., the term
$M_\lm(c)$ below) in the asymptotic of such integrals is given by
the stationary phase method with respect to the parameter $\lm\to \8$ while
the parameter $c$  is {\it fixed} ($c\not=0,\ \pi$). In our case, by a straightforward calculation,   we find out that there are two non-degenerate
critical points of the phase at $t=0$ and $t=\pi/2$. Hence the main
term is a sum of two terms (see equation \eqref{int-K-lm}). We estimate
the remainder {\it uniformly} in $c$ for
$c\not=0,\ \pi$. This is done by reducing the problem to the standard
Beta type integrals. We explain this reduction in Section
\ref{beta-subsect}.

Proposition \ref{rem-hyp-prop} implies that integral
\eqref{int-red2-B} has the following uniform asymptotic expansion in $\lm\in i\br$, $|\lm|\geq \mathcal{S}$ and $c$ ($c\not=0,\ \pi$)  for  {\it fixed} $\tau,\ \tau'$,
\begin{eqnarray}\label{m+r}
l'_\lm(c)=e^{i\frac{\pi}{4}}|\lm|^{-\haf}\cdot
M_\lm(c)+r'_{\lm}(c)\ ,
\end{eqnarray}
where the main term $M_\lm(c)$ comes from stationary points of
the phase at $t=0,\ \pi/2$ and is given by
\begin{eqnarray}\label{int-K-lm}
M_\lm(c)=\left|\sin\left(\frac{c}{2}\right)\right|^{-\lm}+
\left|\cos\left(\frac{c}{2}\right)\right|^{-\lm}\ ;
\end{eqnarray} and for $c\not=0,\ \pi$, the remainder $r'_{\lm}(c)$ satisfies the bound
\begin{eqnarray}\label{rem}
|r'_{\lm}(c)|\leq C|\lm|^{-3/2}|\ln(|\sin(c/2)\cos(c/2)|)|
\end{eqnarray} with a constant $C>0$ depending on $\tau,\
 \tau'$, but not on $c$ and $\lm$.

Let  $m_\lm(c)=|\sin(c/2)|^{-\hlm}|\cos(c/2)|^\hlm$.
After elementary manipulations with \eqref{int-K-lm}, we arrive at
\begin{eqnarray*}\label{int-K-lm-asym}l_\lm(c)=
|\sin(c)|^{{\frac{-\tau-\tau'+\lm-1}{2}}}l'_\lm(c)=&\\
e^{i\frac{\pi}{4}}2^\hlm|\lm|^{-\haf}
|\sin(c)|^{{\frac{-\tau-\tau'-1}{2}}}&\left[m_\lm(c)+m_{-\lm}(c)\right]+|\sin(c)|^{{\frac{-\tau-\tau'-1}{2}}}r'_{\lm}(c). \
\nonumber\end{eqnarray*}  The function $l_\lm$ has the period
equal to $\pi$. We  note that
$m_\lm(c+\pi)=m_{-\lm}(c)$.

In \eqref{def-l-lm} we integrate $l_\lm(c)$
against a function $u$ with a period equal to $\pi$. Hence
we obtain the asymptotic formula \eqref{H-lm-asym}.
\qed

\subsection{Proof of Proposition \ref{G-to-Airy}}
The functional $G_{\lm,n}$ was defined  in \eqref{G-def} through the kernel $k_\lm$ as in \eqref{H-lm-asym}
\begin{eqnarray}\label{I-n-lm}\hskip 1cm
G_{\lm,n}(\phi)=\int_{\br/\pi\bz}\phi(c)|\sin(c)|^{{\frac{-\tau-\tau'-1}{2}}}|\sin(c/2)|^{-\hlm}
|\cos(c/2)|^\hlm\
e^{inc}\ dc
\end{eqnarray}
for $\phi\in C^\8_{even}(S^1)$, $\lm=it\in i\br$, $t\geq \mathcal{S}$, and all $n\in 2\bz_+$.
We consider this integral as a functional on the space of functions $\phi \in C^\8(S^1)$. This functional depends on \lq\lq large" parameters  $\lm$ and $n$, and on axillary parameters $\tau$ and $\tau'$. Our goal is to find a good approximation for values of this functional and give an estimate of the error term.

Let us denote by $S_{\lm, n}(c)= \frac{\lm}{2} (-\ln(|\sin(c/2))| +\ln(|\cos(c/2)|))+
i n c$ the phase of the oscillating
integral \eqref{I-n-lm}  and by $a(c) =
|\sin(c)|^{{\frac{-\tau-\tau'-1}{2}}}$ its amplitude.
Then the functional \eqref{I-n-lm} takes the form
\begin{eqnarray}\label{int2}
G_{\lm,n}(\phi)= \int_{\br/\pi\bz} \phi(c) a(c) e^{S_{\lm,n}(c)} dc \ .
\end{eqnarray}

A direct computation shows that the critical points of the phase
function $S_{\lm,n}$ are solutions of the equation $ \sin(c)= \dl
$, where $\dl = 2in /\lm = 2n/t$.
This shows that the functional \eqref{I-n-lm} has different asymptotic behavior for different values of parameter $\dl$.
 Let us list what we can expect; note that we consider only the case $\dl\geq 0$ (i.e., that $n\geq 0$ and $t\geq\mathcal{S}$).

\begin{enumerate}
  \item For $\dl <1$ the
phase function $S_{\lm, n}$ has two critical points of Morse type;
in this case we can estimate the integral using the stationary phase
method.
  \item When $\dl$ approaches $1$ these critical points collide at the
point $c_0 = \pi/2$. In order to get uniform bounds in this region
we use properties of the Airy function.
  \item When $\dl >1$ the critical points disappear. In this case we
will show that the integral \eqref{I-n-lm} is rapidly decaying.
\end{enumerate}

Our goal is to show that the functional $G_{\lm,n}(\phi)$ can be
 approximated by a functional proportional to the delta function at $c_0$ (i.e., by $A(\lm,n) \phi(c_0)$). We
 will also give explicit uniform bounds for the error term
 $RG_{\lm,n}(\phi) = G_{\lm,n}(\phi) - A(\lm,n) \phi(c_0)$.

We rewrite the phase function $S_{\lm,n}$ in the form $S_{\lm,n}(c)=\hlm S_{\dl}(c)$, where
$\dl=2in/\lm$. We will think about integrals $G_{\lm,n}(\phi)$ as a oscillatory integrals with \lq\lq large" parameter $\lm$ and additional parameter $\dl$.

Using the partition of unity we see that to prove the proposition it is enough to consider
separately two cases:
\begin{enumerate}
  \item The function $\phi$ is supported in a small neighborhood of
the point $c_0=\pi/2$.
  \item The function $\phi$  vanishes in a neighborhood of the point
$c_0 =\pi/2$.
\end{enumerate}

{\it Case 1.} Let $\phi$ be supported in a small enough neighborhood of
the point $c_0=\pi/2$. We claim that for such $\phi$, the following bound holds
\begin{eqnarray}\label{bddd}|G_{it,n}(\phi)-A(it,\dl) \phi(c_0)| \leq  C ||\phi||_{C^2}\cdot  t ^{-\frac{2}{3}}\ . \end{eqnarray}
Here $A(it,\dl) = t^{-\frac{1}{3}}\mathbb{A}(t^{\frac{2}{3}}(\dl-1))=t^{-\frac{1}{3}}\mathbb{A}(t^{-\frac{1}{3}}(2n -t))$, and $\mathbb A$ is  the classical Airy function.

The condition  $1+\eps\geq \dl \geq1-\eps $ implies that $n\asymp |\lm|$. Hence the above bound implies that Proposition  \ref{G-to-Airy} holds for such  $\phi$.

We now specify the size of the support of $\phi$ and prove bound \eqref{bddd}. For any
$0.01>\eps>0$, there exists a
neighborhood $U_\eps\subset [c_0-0.1,c_0+0.1]$ of the point $c_0$ which {\it does not} contain critical
points of $S_\dl$ for $\dl\not \in [1-\eps, 1+\eps]$.  We assume that $\phi$ is supported in this neighborhood for $\eps$ to be specified latter.
Integration by part implies then that for $\dl\not \in [1-\eps, 1+\eps]$, the bound
$|G_{\lm,n}(\phi)|\ll |\lm|^{-N}$ holds for any $N>0$.
Hence we only need to consider the case $1+\eps\geq \dl \geq1-\eps $. We claim that in this case there exists
a change of variables which transforms the integral $G_{\lm,n}(\phi)$ to the Airy type integral. Namely, a direct computation shows that $\frac{\partial}{ \partial c} S_\dl|_{c_0}=\frac{\partial^2}{ \partial c^2} S_\dl|_{c_0}=0$ and $\frac{\partial^3}{ \partial c^3} S_\dl|_{c_0}\not=0$. (In fact,  it is  easy to see that the dependence of $S_\dl$ on $\dl$ is non-degenerate. Namely, the family of functions $\{S_\dl\}$  is a versal deformation of the function $(c-c_0)^3$ in the sense of \cite {Ar}.) We now can quote a classical result on oscillating integrals of the Airy type.  Namely, Theorem 7.7.18, \cite{He} evidently implies the following claim

\begin{claim}{}Let $f\in C^\8(\br^{2})$ be a real valued smooth compactly supported function such that
 $\frac{\partial f}{ \partial x} =\frac{\partial^2f}{ \partial x^2} =0$ and $\frac{\partial^3f}{ \partial x^3} \not=0$
at the point $(x,y)=(0,0)$. Then there exist  $\eps>0$ and smooth real valued functions $a(y)$, $b(y)$ defined on the interval $(-\eps, \eps)$,  such that $a(0)=0,\ b(0)=f(0)$ and
\begin{eqnarray*}\left|\int u(x)e^{i\omega f(x,y)}dx-
e^{i\omega b(y)}\cdot \mathbb{A}\left(a(y)\omega ^\frac{2}{3}\right)\omega ^{-\frac{1}{3}}\cdot u(0)\right|
\leq  C ||u||_{C^2}\cdot \omega ^{-\frac{2}{3}}\ , \end{eqnarray*}
for all real $\omega\geq 1$. Here $\mathbb A$ is  the classical Airy function.
\end{claim}

The above claim  implies  bound \eqref{bddd} for $G_{\lm,n}(\phi)$. Namely,  fix $\eps>0$ such that the above claim is applicable to
$f(x,y)=S_{y+1}(x)$  for $y\in[-\eps,\eps]$ (i.e., for $\dl\in [1-\eps,1+ \eps] $). Let  $U_\eps$ be a  neighborhood of the point $c_0$ which does not contain critical points of $S_\dl$ for $\dl\not \in [1-\eps, 1+\eps]$. We assume that $supp(\phi)\subset U_\eps$. Applying the above claim for $\dl\in [1-\eps, 1+\eps]$, we obtain the bound \eqref{bddd}.

{\it Case 2.} Let $\phi$ be a function vanishing in a neighborhood of the point
$c_0 =\pi/2$. In this case we have  upper
bounds
\begin{eqnarray}
|G_{\lm,n}(\phi)|\leq\left\{
 \begin{array}{ll}
 C'||\phi||_{C^2}\cdot|\lm|^{-\haf}, & \hbox{for}\ \dl>0.9, \\
 C'_{N}||\phi||_{C^N}\cdot|\lm|^{-N}, & \hbox{for}\ 0<\dl\leq 0.9,
 \end{array}
\right.
\end{eqnarray} for any $N>0$ and some constants $C',\ C'_{ N}$, which could
be explicitly bounded in terms of $\tau$ and $\tau'$. These bounds
immediately follow from the van der Corput lemma and integration by parts as explained in Section \ref{intgr-parts-end}.
\qed


\appendix
\section{Beta integrals}\label{red-hyp-int-sect}

In this appendix we explain how to prove asymptotic expansion for certain oscillating integrals which we call Beta integrals. We use these asymptotic in the proof of  Lemma~\ref{claim-H-lm}.

\subsection{Beta integrals}\label{beta-subsect}
Fix a function  $h\in C^\8(\br)$ such $h(0)=0$,  $h'>0$. Fix $\s,\ \s'\in\bc$ such that  $Re(\s),\ Re(\s')>-1$ and $Re(\s)+Re(\s')=-1$. (In fact, in this paper we will need only the case  $Re(\s)=Re(\s')=-\haf$). We consider following integrals

\begin{eqnarray}\label{par-1} \mathrm{H}_{\lm,c}(\phi)=
\int_\br |h(t-c)|^{\s+\lm}|h(t+c)|^{\s'+\lm}\phi(t)
dt\ ,
\end{eqnarray} where $\phi\in C^\8(\br)$,  and $\lm\in i\br$.  We are
interested in the {\it uniform} asymptotic of such integrals in $c$, $c\not=0$, and   for
$|\lm|$ sufficiently large. Moreover, we will assume that both $supp(\phi)$ (containing $0$)
and values of $c$ are sufficiently small, depending on the function $h$.

We write the integral $\mathrm{H}_{\lm,c}(\phi)=\int _\br  \phi(t) a_{\s,\s'}(c;t) e^{\lm S(c;t)} dt $ in the standard form customary in the stationary phase method. Here $\phi(t) a_{\s,\s'}(c;t)$ is the amplitude and $S(c;t)$ is the phase in this oscillating integral, both depending on the parameter $c$ and some auxiliary parameters $\s,\ \s'$ which we consider fixed. For any {\it fixed} $c\not=0$ and smooth $\phi$ of compact support, one can obtain the asymptotic in
$|\lm|\to\8$ for $\mathrm{H}_{\lm,c}(\phi)$ from the stationary phase method (see \cite[Theorem 7.7.6]{He}).
We choose the range of the parameter $c$ and the support of $\phi$ small enough so that  for all $c\not=0$, the following conditions are satisfied. There exists the unique  critical point $t_{c}$ (in variable $t$) of the phase  $S(c;t)$,  this critical point is non-degenerate, and it is disjoint from singularities of the amplitude
$a_{\s,\s'}$ at points $t=\pm c$ (in fact if $h$ is odd, as in our case, then $t_c=0$ for all $c\not=0$). We denote by
$\mathrm{H}^0_{\lm,c}(\phi)$ the main term of the contribution from the critical point $t_c$  to the asymptotic of
$\mathrm{H}_{\lm,c}(\phi)$ given by the stationary phase method. (In particular, we will show that for large $|\lm|$ and  fixed  $c$, $|\mathrm{H}^0_{\lm,c}(\phi)|=A|\lm|^{-\haf}$ and $|\mathrm{H}_{\lm,c}(\phi)-\mathrm{H}^0_{\lm,c}(\phi)|\leq B|\lm|^{-3/2}$.)

Our aim is to obtain a meaningful  bound for the remainder
\begin{eqnarray*}R\mathrm{H}_{\lm,c}(\phi)=\mathrm{H}_{\lm,c}(\phi)-\mathrm{H}^0_{\lm,c}(\phi)\ , \end{eqnarray*}which is {\it uniform} in $\lm$ and $c$.
Recall that we set $\mathcal{S}=2(|\tau|+|\tau'|)+1$. We claim the following bound

\begin{prop}{rem-hyp-prop}Fix $h\in C^\8(\br)$ as before. There are constants $C_1,\ C_2>0$, and intervals  $(-\epsilon ,\ \epsilon)$ and $[-d,\ d]$ depending on the function $h$, such that the remainder satisfies the bound
\begin{eqnarray}\label{hyper-rem-app}
|R\mathrm{H}_{\lm,c}(\phi)|\leq C_1||\phi||_{C^{1}}\cdot|\lm|^{-\frac{3}{2}}+C_2||\phi||_{C^{2}}\cdot|
\ln|c||\cdot|\lm|^{-2}\
\end{eqnarray}for any  $|\lm|\geq \mathcal{S}$, $c\in (-\epsilon,\ \epsilon),\ c\not=0$, and for any smooth function $\phi$ such that $supp(\phi)\subset [-d,\ d]$.
\end{prop}

In fact the method we present allows one to give the asymptotic expansion to any order with the explicit bound on the remainder.

\subsection{Proof of Proposition \ref{rem-hyp-prop}}
We show that it is enough to consider the special case of $h(t)=t$.
Namely, we claim
there exists a smooth change of variables $(t, c)$ to the new set
of variables $(x,a)$, where $c$ depends on $a$ only, such that
it transforms the kernel function
$|h(t-c)|^{\s+\lm}|h(t+c)|^{\s'+\lm}$ to the homogenous kernel
$|x-a|^{\s+\lm}|x+a|^{\s'+\lm}$ times some smooth function mildly
depending on $a$.

Let $g\in C^\8_c(\br)$ be a function such that $h(t)=tg(t)$ and
$g(0)>0$. We denote by $f(t,c)=h(t-c)h(t+c)$. The necessary change
of variables is given by the following lemma.

\begin{lem}{change-hyp}There exists a change of variables
$(x,a)=(x(t,c),a(t,c))$ in a neighborhood of the point $(0,0)$
such that
\begin{enumerate}
 \item The variable $a$ is a function of $c$ only,
 \item $f(t,c)=(x+a)(x-a)$ in new coordinates, and
 \item $h(t-c)=(x-a)g_1(x,a)$ and $h(t+c)=(x+a)g_2(x,a)$,
 where $g_1$ and $g_2$ are smooth functions not vanishing
 near the point $(0,0)$.
\end{enumerate}
\end{lem}

Using this lemma,  we can rewrite the integral
\begin{eqnarray}\label{100} \mathrm{H}_{\lm,c}(\phi)=
\int_\br |h(t-c)|^{\s+\lm}|h(t+c)|^{\s'+\lm}\phi(t)
dt\ = \\ \int_\br
|x-a|^{\s+\lm}|x+a|^{\s'+\lm} \psi(x) dx\ ,\nonumber
\end{eqnarray} where $\psi$ is a smooth function such that $\psi(0)=\phi(0)$ and $C^n$-norms of $\psi$ are bounded by those of $\phi$. Explicitly   $\psi(x)=\phi(t(x,a))|g_1(x,a)|^\s|g_2(x,a)|^{\s'}\left|\frac{\partial x}{\partial t}\right|$.

We introduce integrals
\begin{eqnarray}\label{231}
H_{\lm,a}(\psi)=\int_\br|x-a|^{\s+\lm}|x+a|^{\s'+\lm} \psi(x) dx\ .
\end{eqnarray}

Lemma \ref{change-hyp} implies  that
$ \mathrm{H}_{\lm,c}(\phi)= H_{\lm,a}(\psi)$ for an appropriate function $\psi$ (see \eqref{100}). Here parameters $c$ and $a$ are related via the change of variables in Lemma~\ref{change-hyp}.

The integral $H_{\lm,a}(\psi)$ also has an asymptotic
expansion (in $\lm$ for every {\it fixed} $a$) with the main term
$H^0_{\lm,a}(\psi)$ given by the stationary phase method at $x=0$, and a
remainder $RH_{\lm,a}(\psi)$. We want to compare asymptotic
expansions of $\mathrm{H}_{\lm,c}(\phi)$ and of $H_{\lm,a}(\psi)$.
Our considerations are based on the well-known {\it invariancy} of
terms obtained by the stationary phase method (see \cite{Ar},
\cite{St}). Namely, we have
$H^0_{\lm,a}(\psi)=\mathrm{H}^0_{\lm,c}(\phi)$. Since integrals
themselves are also equal we have the equality of remainders
$RH_{\lm,a}(\psi)= R\mathrm{H}_{\lm,c}(\phi)$. Hence, we can use the
estimate for the remainder for the  integral $H_{\lm,a}$ which we obtained in \eqref{hyper-rem-app3}, Corollary  \ref{rem-hyp-cor}.

Parameters  $a$ and $c$ belong to a bounded set. Hence $C^N$-norms of  $\psi$ could be
bounded independently of $a$ in terms of $||\phi||_{C^{N}}$ and of $||h||_{C^{N}}$. This implies that the constant in the bound \eqref{hyper-rem-app} for the remainder $R{H}_{\lm,a}(\psi)$ could be chosen  independently of $c$. Hence we reduced the proof of bound \eqref{hyper-rem-app} for general function $h$ to the special case $h(t)=t$. This special case is dealt with in the next section (see Corollary  \ref{rem-hyp-cor}). This finishes the proof of Proposition \ref{rem-hyp-prop}. \qed

\subsection{Standard Beta integrals}\label{stand-beta-sect}  Consider following standard Beta integrals
\begin{eqnarray}\label{Hlm} H_{\lm,\s,\s'}(\phi)=\int_\br
|y-1|^{\s+\lm}|y+1|^{\s'+\lm} \phi(y) dy\ ,
\end{eqnarray}
where   $\phi\in C^\8(\br)$, $\lm\in i\br$, and $\s,\ \s'$ are as before.  We apply the stationary phase method and the elementary method of integration by parts as described in Section \ref{intbyparts} in order to obtain the following bound.

Let $R=\br\setminus [-0.5,\ 0.5]$ and $\xi=y\frac{\partial}{\partial y}$. The phase function in integral \eqref{Hlm} has the unique stationary point at $y=0$ which is non-degenerate. Let $H^0_{\lm,\s,\s'}(\phi)$ be the main term in the asymptotic of $H_{\lm,\s,\s'}(\phi)$ as $|\lm|\to\8$ (i.e., $H^0_{\lm,\s,\s'}(\phi)= \al\phi(0)\cdot|\lm|^{-\haf}$ with $\al=\left(\frac{\pi}{i}\right)^\haf$ given by the stationary phase method).

\begin{lem}{rem-hyp-lem}There are constants $ C_1,C_2>0$ such that the bound
\begin{eqnarray*}
\left|H_{\lm,\s,\s'}(\phi)-H^0_{\lm,\s,\s'}(\phi)\right|\leq C_1||\phi||_{C^{1}([-0.9, 0.9])}\cdot|\lm|^{-\frac{3}{2}}+C_2 RH(\phi)\cdot |\lm|^{-2}
\end{eqnarray*} holds for any  $|\lm|\geq \mathcal{S}$,  and for any smooth compactly supported function $\phi$.
Here the reminder is given by $RH(\phi)=\int\limits_R\sum\limits_{i=0}^{2}|\xi^i(\phi)|\frac{dy}{|y|}$.
\end{lem}

\Proof It is enough to treat separately the case of $\phi$ supported near zero (e.g., in the interval  $[-0.9, 0.9]$)  and that of $\phi$ vanishing near zero (e.g., vanishing on $ [-0.5,\ 0.5]$).

{\it Case 1.} Function  $\phi$ supported near zero. The  stationary phase method (see \cite[Theorem 7.7.6]{He})  implies that
\begin{eqnarray}\label{stat-phase} |H_{\lm,\s,\s'}(\phi)-H^0_{\lm,\s,\s'}(\phi)|\leq C_1||\phi||_{C^1}\cdot |\lm|^{-\frac{3}{2}}\ ,
\end{eqnarray}
with an explicit constant $C_1$. Such a bound is enough for our purposes.

{\it Case 2.} Function  $\phi$ vanishes near zero. We rewrite  the integral $H_{\lm,\s,\s'}(\phi)$  in the form $I_F$  from \eqref{IF}, Appendix \ref{corput}, with
\begin{eqnarray}\label{Ilm} F(y;\lm,\s,\s')=y|y-1|^{\s}|y+1|^{\s'}|y-1|^{\lm}|y+1|^{\lm} ,
\end{eqnarray}
and the form $\omega=dy/y$.

Consider the vector field $\xi=y\frac{\partial}{\partial y}$. A straightforward computation shows that
 $G:=\xi(F)/F=\lm(\frac{y}{y+1}+\frac{y}{y-1})+g_{\s,\s'}(y)$, where the function
$g_{\s,\s'}$ is bounded on the set $R=\br\setminus [-0.5,\ 0.5]$. Hence, for $|\lm|\geq \mathcal{S}$, the function $H=G\inv$ is uniformly bounded in $\lm$ and $y\in \br\setminus [-0.5,\ 0.5]$. Moreover, if we make a change of variable $z=y\inv$, then the function $H$ and the vector field $\xi$ are smooth on the interval $J=[-1,\ 1]$ (including at zero, after extending $H$ and $\xi$ by continuity).  Via compactness, this implies that all functions $\xi^i(H)$ are uniformly bounded (in the coordinate $z$) on $J$, and hence are bounded on $\br\setminus [-0.5,\ 0.5]$ (in the original coordinate $y$). This allows us to estimate the integral $I_F(\phi)$ and finishes the proof of the lemma.
\Endproof

We will use the bound described in the lemma in order to estimate the integral $H_{\lm,a}$ as defined in \eqref{231}. Clearly we can reduce the integral $H_{\lm,a}$ to the standard Beta integral  $H_{\lm,\s,\s'}$. Namely,
\begin{eqnarray*}
H_{\lm,a}(\psi)=\int_\br|x-a|^{\s+\lm}|x+a|^{\s'+\lm} \psi(x) dx=\\ |a|^{\s+\s'-1+2\lm}\int_\br |y-1|^{\s+\lm}|y+1|^{\s'+\lm} \psi(ay) dy\ .
\end{eqnarray*}

Let $H^0_{\lm,a}(\psi)$ be the main term in the asymptotic  of  $H_{\lm,a}(\psi)$ which is given by the stationary phase method for $a\not=0$ fixed. Applying the above lemma to the last integral we obtain the following bound.

\begin{cor}{rem-hyp-cor}Let $\psi$ be a compactly supported smooth  function. There are constants $C_3,C_4>0$, depending on $\psi$ such that the bound
\begin{eqnarray}
\label{hyper-rem-app3}\hskip 1cm
\left|H_{\lm,a}(\psi)-H^0_{\lm,a}(\psi)\right|\leq C_3||\psi||_{C^{1}([-0.9, 0.9])}\cdot|\lm|^{-\frac{3}{2}}+C_4|\ln(a)|\cdot|\lm|^{-2}\
\end{eqnarray} holds for all   $|\lm|\geq \mathcal{S}$ and $a\in(0,\ 0.1]$.
\end{cor}

We have $H^0_{\lm,a}(\psi)=|a|^{\s+\s'-1+2\lm}\al|\lm|^{-\haf}\psi(0)$. Note that we assumed that $Re(\s+\s'-1+2\lm)=0$ and hence $|H^0_{\lm,a}(\psi)|=|\al\psi(0)|\cdot|\lm|^{-\haf}$.

\Proof Let $supp(\psi)\subset[-A,\ A]$ and denote by $\psi^a(y)=\psi(ay)$. We note that  $\sup|\xi^i(\psi^a)|\leq\sup |\xi^i(\psi)|$ for any $a\in(0,\ 0.1]$. Hence we have
\begin{eqnarray*}\label{bound} |RH(\psi^a)|&\leq& C_1|\lm|^{-n}\sum_i\int|F|\xi^i(\psi(ay)u(y))||\omega| \\ &\leq& C_2|\lm|^{-n}\sum_i\int_\haf^{a\inv A}|F||\omega|\leq C_3|\lm|^{-n}|\ln(a)|\ ,
\end{eqnarray*} for any $n$ and for some explicit constants $C_{1,2,3}$ depending on derivatives of $\psi$. Here we use the fact that $|F|$ is bounded as $y\to\pm\8$ and that $\omega=dy/y$.
\Endproof

\subsubsection{Proof of Lemma \ref{change-hyp}}\label{var-chng-hyper}
The proof is based on the theory of normal forms of differentiable
functions and
on  Hadamard's lemma (see \cite{Ar}, \cite{Ma}).

Consider a smooth family of functions $f(t,c)=h(t-c)h(t+c)=(t^2-c^2)g(t-c)g(t+c)$,
 where we view $t$ as a variable and $c$ as a parameter. For $c=0$ the function
 $f(t,0)$ is equivalent (under a smooth change of variable $t$) to the function $t^2$.
 The theory of versal deformations then implies that there is a change of variable
 $x=x(t,c)$ such that $f(x,c)=u(c)+x^2$ for some smooth function $u$ (see \cite{Ar}). On the other
 hand, the differential of $f(t,c)$ with respect to $c$ vanishes for all
 $t$ and $f(0,c)<0$. This implies that we can write $u(c)=-c^2\tilde u^2(c)$
 with $\tilde u(c)>0$. Hence there exists a new parameter $a=a(c)$ such that
 $f(x,a)=x^2-a^2=(x-a)(x+a)$.

By Hadamard's lemma (see  \cite{Ma}) $h(t-c)$ is divisible by $(x-a)$ since these functions have the
same zeroes (one of the branches of zero set for the function
$f(x,a)=x^2-a^2$). Hence we can write $h(t-c) = (x-a) g_1(x,a)$.
It is clear that $g_1$ is invertible near $0$. Similarly for the
function $h(t+c)$.\qed


\section{Integration by parts and Van der Corput lemma}\label{corput}
\subsection{Integration by parts}\label{intbyparts} We want to study integrals of the form
\begin{eqnarray}\label{IF}I_F(\phi)=\int_\br F(y;\lm,r)\phi(y)\omega\ ,
\end{eqnarray}
where $\omega$ is a one-form in $y$, $F$ is a certain kernel depending on a large parameter $\lm\in \bc$ and on some additional (multi)parameter $r\in \br^m$, and $\phi$ is of compact support. We would like to obtain estimates of $I_F$ for $|\lm|\to\8$. We are interested in uniform in $r$ estimates given in terms of $C^k$-norms of the function $\phi$ (i.e., we want to estimate a $C^k$-norm of the functional $I_F$).
We have the following elementary method based on the integration by parts.

First we note that there is a trivial estimate for integral \eqref{IF} by the absolute value: $|I_F(\phi)|\leq R_F(\phi)$, where $R_F(\phi)= \int_\br |F(y;\lm,r)\phi(y)||\omega|$.
We use the integration by parts to bootstrap this estimate.

Suppose we are given a vector field $\xi$ on $\br$ and a function $H=H(y,\lm, r)$ such that
\begin{description}
  \item[\it (i)] $H\cdot\xi( F)=\lm \cdot F$,
  \item[\it (ii)]  $H$ is a smooth in all variables, and for some $n>0$,  absolute values of functions $H,\ \xi H,\dots, \xi^n H$ are  bounded by a constant $C>0$,  uniformly in all parameters,
  \item[\it (iii)] $\xi \omega=0$.
\end{description}

\begin{prop}{intbyparts-prop} For $\xi$ and $H$ as above, we have the following bound
\begin{eqnarray}\label{prop-intbparts-bnd} |I_F(\phi)|\leq |\lm|^{-n}\cdot C^n\sum\limits_{i=0}^{n}R_F(\xi^i \phi)\ .
\end{eqnarray}
\end{prop}
\Proof
 We have the following functional equation
\begin{eqnarray} I_F(\phi)=-\lm\inv\cdot I_F(\xi(H\phi))\ .
\end{eqnarray}
Indeed, we have
\begin{eqnarray*} I_F(\xi(H\phi))=\int F\cdot \xi(H\phi)\omega=-\int\xi(F)H\phi\ \omega=-\lm\int F\phi\omega=-\lm I_F(\phi)\ .
\end{eqnarray*}
Iterating this we obtain $I_F(\phi)=(-1)^n|\lm|^nI_F(D^n(\phi))$, where $D(\phi)=\xi(H\phi)$. Clearly we have
$D^n(\phi)=\sum\limits_{0\leq i_0,\dots,i_{n+1}\leq n}\left[H^{i_0}\cdot(\xi(H))^{i_1}\cdot(\xi^2(H))^{i_2}\dots(\xi^n(H))^{i_{n}}\right]\cdot \xi^{i_{n+1}}( \phi)\ ,$
where the summation is over an appropriate set of indexes. Hence we arrive at the desired bound
\begin{eqnarray} |I_F(\phi)|\leq |\lm|^{-n}\cdot C^n\sum\limits_{i=0}^{n}\int |F||\xi^i (\phi)||\omega|\ .
\end{eqnarray}
\Endproof

\subsection{Van der Corput lemma}

 Let $f$ be a real valued smooth function on the
interval $[a,b]$, and $F(x)=e^{if(x)}$.   Consider the following integral
\begin{equation}\label{corput-int}
 I(f,\phi):=I_F(\phi)= \int_a^b e^{i f(x)}\phi(x)dx\ .
\end{equation}

The bound \eqref{prop-intbparts-bnd} implies the following

\begin{cor}{}Let $f=t\al$, where $t>1$ is a real parameter and $\al$ is a smooth function such
that $\al'$ has no zeroes on the support $supp(\phi)\Subset(a,b)$ of $\phi$. Then the following bound holds
\begin{equation}\label{int-parts}
 |I(t\al, \phi)|\leq C_Nt^{-N}
\end{equation}
for any $N>0$, and a constant $C_N$ depending on $\al$ and $\phi$.
\end{cor}

Let $I(f,\phi)$ be as in \eqref{corput-int}. Consider  the case when $f'$ has zeroes.
For an integer $k\geq 1$ denote by
$m_k(f)=\min\limits_{x\in[a,b]}|f^{(k)}(x)|$ and let $M(\phi)=
|\phi(b)|+\int_a^b|\phi'(x)|dx$ be the variance of $\phi$. We have the
following general estimate
essentially due to van der Corput (see \cite[p. 332]{St}).

\begin{lem}{corput-l}Let $k\geq 1$ be such that $m_k(f)>0$.
There exists a constant $c_k$ such that the following bound holds
\begin{equation*}\label{corput-bd}
 \left|I(f,\phi)\right|\leq c_k\cdot
 m_k(f)^{-\frac{1}{k}}\cdot M(\phi)
\end{equation*}
provided
\begin{enumerate}
 \item $k\geq 2$, or
 \item $k=1$ and $f'$ is monotone on $[a,b]$.
\end{enumerate}
The constant $c_k$ depends only on $k$ and is independent of $\phi$,
$f$ and of the interval
$[a,b]$.
\end{lem}

We use this lemma with $k=1$ or $2$, so we can assume that $c_k$
is a universal constant.

\subsection{}\label{intgr-parts-end} Throughout the paper we consider integrals of the form $\int u(x)|x|^{-it}e^{i s\cdot g (x)}dx$.
In this section we explain how to obtain meaningful upper bounds for
these integrals. We claim
that the necessary type of bounds follow directly from  the integration by parts and from the van der
Corput lemma.

Let
\begin{equation}\label{corp-1}
 I( s ,t )=\int_{-1}^1u(x)|x|^{-\haf-it }e^{i s\cdot  g (x)}dx\ ,
\end{equation}
where we assume that $1\leq t \leq s $, $g $ is smooth and
monotonic, $0.99<g '(x)<1.01$ (i.e., bounded away from $0$ and
$\8$), $|g ^{''}(x)|\leq \haf$ for all $x$ (this insures that
there is no degenerate critical points of the phase), and $u$ is
smooth of compact support in $(-1,1)$.

There is a simple bound if the phase has no critical points. Let us
denote by $b$ the ratio $b= s /t $. Integration by parts shows that if
 the phase
function $bg (x)-\ln|x|$ in the integral \eqref{corp-1} has no
critical points (e.g., $|t |\gg| s |$) then the bound \eqref{int-parts} reads as
\begin{equation}\label{corp-N}|I( s ,t )|\leq C_N|t |^{-N}\end{equation}
 for any $N>0$ and some constant $C_N$ depending on $N,\ u$ and $g $.

In the complementary situation we have

\begin{lem}{}Under the above assumptions on $g $, the following uniform
 bound holds
\begin{equation*}\label{corp-2}
 |I( s ,t )|\leq B  s ^{-\haf}\ ,
\end{equation*} where the constant $B$ is independent of $ s $ and of
$t $.
\end{lem}

\Proof We denote by $a$ the ratio $a=t / s $ and consider
the integral over the interval $(0,1)$ (and the similar integral
over $(-1,0)$)
\begin{equation*}\label{corp-11}
\int_{0}^1u(x)|x|^{-\haf}e^{i s  (g (x)-a\ln|x|)}dx\ .
\end{equation*} We are interested in the uniform (in $ s $) bound
for this integral for the values of the parameter $a$ satisfying
the bound $ s \inv\leq a\leq 1$.

In order to apply the van der Corput lemma, we break the interval
$(0,1)$ into $4$ intervals $J_1=(2a, 1)$, $J_2=(\haf a, 2a)$,
$J_3=(\haf s \inv,\haf a)$ and $J_4=(0,\haf s \inv)$ (for $a\geq
\haf$ the first interval is missing). Denote by
$f_a(x)=g (x)-a\ln|x|$, $\phi(x)=u(x)|x|^{-\haf}$ and consider
the corresponding
integrals $I_j( s , a)=\int_{J_j}u(x)|x|^{-\haf}e^{i s
f_a(x)}dx$.

On the interval $J_1$ we have $| s  f_a'(x)|\geq  s $. Hence from the
van der Corput lemma (with $k=1$) we have $|I_1( s ,a)|\leq
B_1 s ^{-1}$.

On the interval $J_2$ the phase $f_a$ has zero of the first
derivative, but satisfies the bound $| s  f_a^{''}(x)|> \haf
a\inv s $ and $M(\phi)\leq 10|a|^{-\haf}$. Hence on the interval $J_2$
the van der Corput lemma with
$k=2$ implies $|I_2( s ,a)|\leq B_2 s ^{-\haf}$.

To bound the integral $I_3( s ,a)$, we note that $| s  f_a'(x)|\geq
\haf s $ and that the variation of the amplitude satisfies $M(\phi)\leq|\haf
a|^{-\haf}+\int^{\haf a}_{\haf s \inv}|x|^{-3/2}dx\leq c  s ^{\haf}$ on $J_3$. The
van der Corput lemma with $k=1$ implies that $|I_3( s ,a)|\leq B_3
 s ^{-\haf}$.

Bounding the integral over $J_4$ by the integral of the absolute
value, we see that trivially $|I_4( s ,a)|\leq B_4 s ^{-\haf}$ .
\Endproof


\end{document}